\documentclass{article}

\usepackage{amssymb,latexsym,amsmath, amscd}

\newtheorem{theorem}{Theorem}[section]
\newtheorem{corollary}[theorem]{Corollary}
\newtheorem{lemma}[theorem]{Lemma}
\newtheorem{proposition}[theorem]{Proposition}
\newtheorem{definition}[theorem]{Definition}

\newtheorem{remark}[theorem]{Remark}

\begin{document}

\renewcommand{\contentsname}{ }

\date{This version:  30.09.2009}

\title{A characterization of sub-riemannian spaces
 as length dilatation structures constructed via coherent projections}

\author{Marius Buliga \\
\\
Institute of Mathematics, Romanian Academy \\
P.O. BOX 1-764, RO 014700\\
Bucure\c sti, Romania\\
{\footnotesize Marius.Buliga@imar.ro}}

\date{This version:  10.10.2009}

\maketitle

\begin{abstract}
We introduce length dilatation structures on metric spaces, tempered dilatation
structures and coherent projections and explore the
relations between these objects and the Radon-Nikodym property and
 Gamma-convergence of length functionals. Then we show that 
the main properties of sub-riemannian spaces can be obtained from pairs of 
length dilatation structures, the first being a tempered one and the second 
obtained via a coherent projection.  Thus we get an intrinsic, synthetic, 
 axiomatic description of sub-riemannian geometry, which transforms 
  the classical construction of a Carnot-Carath\'eodory distance on a 
  regular sub-riemannian manifold into a model for this abstract sub-riemannian 
  geometry. 
\end{abstract}



\tableofcontents


\section{Introduction}

Sub-riemannian geometry is the study of non-holonomic spaces (introduced by Vr\u anceanu \cite{vra1}, \cite{vra2} in 1926) 
endowed with a Carnot-Carath\'eodory distance. Such spaces appear in 
applications to thermodynamics (the name "Carnot-Carath\'eodory 
distance" is inspired by the  work of Carath\'eodory \cite{cara} (1909) concerning 
a mathematical approach to Carnot work in thermodynamics), in 
  non-holonomic dynamics (see the survey Vershik and Gershkovich \cite{vershik1}),  
in the study of hypo-elliptic operators  H\"ormander \cite{hormander}, in 
harmonic analysis on homogeneous cones  Folland, Stein \cite{folstein}, and
as boundaries of CR-manifolds. 

In several foundational papers on sub-Riemannian geometry,
among them Mitchell \cite{mit}, 
Bella\"{\i}che \cite{bell}, the  paper of Gromov asking for an
 intrinsic point of view for sub-riemannian geometry  \cite{gromovsr}, Margulis,
 Mostow \cite{marmos1}, \cite{marmos2}, dedicated to Rademacher theorem for
 sub-riemannian manifolds and to the construction of a tangent bundle of such
 manifolds, and Vodopyanov \cite{vodopis} \cite{vodopis2}, Vodopyanov and
 Karmanova \cite{vodokar},  fundamental results concerning the intrinsic
 properties of sub-riemannian spaces endowed with the Carnot-Carath\'eodory
 distance  were proved using differential geometry tools, which are  not intrinsic to
sub-Riemannian geometry.

The point of view of Gromov in \cite{gromovsr} is  that the only intrinsic 
object on a sub-riemannian  manifold is the Carnot-Carath\'eodory  distance. 
The underlying differential structure of the 
manifold is  then clearly not intrinsic. Nevertheless, in all proofs in the
before mentioned papers on the fundamentals of sub-riemannian geometry this
differential structure is used in order to prove intrinsic statements. 

Among the fundamental results in sub-riemannian geometry, or more general non-holonomic 
geometry, a particular position occupies the   result concerning  the nilpotent 
group structure of the tangent space to a point  of a regular  non-holonomic manifold 
(for an evolution of this subject see \cite{vershik2}, \cite{agrachev}).  According to the 
introduction of the Agrachev and Marigo paper \cite{agrachev}, non-holonomic tangent functors 
appear appear via constructions involving nilpotent or graded approximations of the 
geometrical objects to be studied.  We cite from \cite{agrachev} page 112, 3rd paragraph: 
"A weak point of these constructions is their heavy dependence on the choice of coordinates. 
Because of that, the approximation looks like an auxiliary technical tool rather than a fundamental functorial operation;  the geometric insight and the application of geometric machinery are highly impeded."  Agrachev and Marigo propose therefore an intrinsic construction of the tangent bundle 
of a non-holonomical manifold.  Their notion of "intrinsic" means "coordinate free", in the frame of differential geometry of manifolds. 

Notice that studies like \cite{vershik2}, \cite{agrachev}, show that in fact a very important geometrical 
object of sub-riemannian geometry does not depend on the Carnot-Carath\'eodory distance in any direct 
way.  

In conclusion, for  sub-riemannian geometry there are two meanings of the word "intrinsic" : 
\begin{enumerate} 
\item[(a)] the Carnot-Carath\' eodory distance as the only intrinsic object permits to formulate 
several important results, but the proof of these results is based on the constructions and 
approximations mentioned by Agrachev and Marigo, 
\item[(b)]  if we see sub-riemannian geometry as a subspecies of non-holonomic geometry, then 
"intrinsic" is only the non-holonomic distribution and differential geometry is accepted as an 
"intrinsic"  tool. 
\end{enumerate}

The interpretation of this situation is that the Carnot-Carath\'eodory distance is necessary, but not sufficient for an axiomatization of sub-riemannian geometry. On the other side differential geometry 
is certainly sufficient to describe what sub-riemannian geometry is, but it seems that non-holonomic 
manifolds are just models of a  non-holonomic geometry.  In few words,   the Carnot-Carath\'eodory distance is not enough and the whole formalism of differential geometry is too much for sub-riemannian geometry.

We tried (\cite{buligadil1}, 
\cite{buligadil2}, \cite{buligasr}) to find an intrinsic  frame, not using differential geometry, 
 in which sub-Riemannian geometry would be a model. The main sources of inspiration were   the last 
section of the paper by Bella\"{\i}che \cite{bell}, the intrinsic point 
of view of Gromov \cite{gromovsr}, the Pansu derivative  \cite{pansu} 
and the author previous attempt to understand the operation 
of differentiation from a topological point of view \cite{buliga0}.

Besides  the Carnot-Carath\'eodory distance, there is another object which persistently appears 
in all studies of sub-riemannian geometry: the (anisotropic) dilatations.

We first proposed the notion of dilatation 
structure, studied in \cite{buligadil1}. A dilatation structure encodes 
the approximate self-similarity of a metric space and it induces 
a metric tangent bundle  with group operations in each fiber 
(tangent space to a point), which make it (the tangent space) into a
conical group.  Conical groups generalize Carnot groups. The affine geometry of conical groups was then studied in 
\cite{buligadil2}. In  
\cite{buligasr}  it is shown that regular sub-riemannian manifolds 
admit  dilatation structures constructed via normal frames. In that  paper 
we tried to  minimize 
the contribution of classical differential calculus in the proof of the
basic results in sub-riemannian geometry, by showing that in fact the 
differential calculus on the underlying differential manifold of the 
 sub-riemannian space is needed only for 
proving that normal frames exist, which implies the existence of 
dilatation structures associated to the Carnot-Carath\'eodory distance.

The point of view of this paper is that sub-riemannian geometry may be described
by a set of axioms concerning dilatation structures. It is true that 
this viewpoint is less general than Gromov's (there are more intrinsic objects 
than the Carnot-Carath\'eodory distance). Nevertheless, in this approach we
renounce at the differential structure (of the manifold) and we replace it 
with something which is much weaker, a dilatation structure.

In \cite{buligaultra} we showed that there are many 
dilatation structures  on ultrametric spaces.  The distance on these metric spaces is not a length 
distance, therefore such dilatation structures are different from the ones appearing in  sub-riemannian geometry. We could then imagine that (a generalization of the) sub-riemannian geometry is the  study of dilatation structures on length metric spaces. 

With  this motivation we propose here the notion of a length dilatation structure (section \ref{seclds})  with the Radon-Nikodym property (RNP) (section \ref{radon}). The dilatation structure of a regular sub-riemannian manifold 
is a length dilatation structure with RNP.  Beside these, are there  any other axioms which we have 
to add in order to obtain a class of dilatation structures which describes ths sub-riemannian geometry?

The answer (theorem \ref{mainsrthm}) is that they are length dilatation structures (definition
\ref{deflds}) and they are constructed with the help of coherent projections (definition \ref{defcoh}) 
and tempered dilatation structures (definition \ref{dtempered}). 

Tempered dilatation structures have the property that  for any point $x$  of the space, the 
dilatation based at $x$ is bi-lipschitz, in a uniform manner with respect to 
the magnification $\varepsilon$ and the base point $x$. This property describes riemannian spaces,  
but not general sub-riemannian spaces. In order to obtain sub-riemannian spaces we also need 
coherent projections, which are objects generalizing non-holonomic distributions.

 The ingredients 
of the classical construction of a sub-riemannian manifold are  
a riemannian manifold and a distribution. From these ingredients a new distance
is constructed: the Carnot-Carath\'eodory or sub-riemannian distance. 
The construction proceeds then further, by showing various convergences of differential geometrical
quantities (vector fields, deformed riemannian metrics) to corresponding 
quantities which give a structure to the metric tangent space at a point from
the space (initial manifold) endowed with the sub-riemannian distance.  
This construction is generalized here to dilatation structures by replacing 
distributions by coherent projections.

Consider $M$ a real smooth $n$-dimensional manifold. We may think 
 in the first instance that instead of a distribution, which is a map associating
 to any point $x \in M$ a subspace $\displaystyle D_{x}  \subset T_{x}M$,  we 
 use a field of projections 
 $$\displaystyle Q^{x}: T_{x}M \rightarrow T_{x}M \quad , \quad Q^{x} \, 
 T_{x}M \, = \, D_{x} \quad , \quad Q^{x}\, Q^{x} \, = \, Q^{x} $$
 But where these projections are coming from and why do we think about them as
 more interesting as distributions? Let us denote by 
 $\displaystyle \bar{\delta}^{x}_{\varepsilon} u \, = \, \varepsilon u$ the usual
 multiplication by positive scalars in the tangent space of $M$ at $x$. Suppose
 that the distribution $D$ is spanned by a family of vector fields which induces 
 by the Chow condition 
 a normal frame $\displaystyle \left\{ X_{i} \mbox{ : } i = 1, ... , n \right\}$, definition \ref{defnormal},
  and a non-isotropic dilatation
$$\delta^{x}_{\varepsilon} \left( \sum_{i=1}^{n} a_{i} X_{i}(x)\right) \  = \
\left( \sum_{i=1}^{n} a_{i} \varepsilon^{deg X_{i}}  X_{i}(x) \right)$$
  as in theorem \ref{structhm}, then 
 $$Q^{x} u \ = \  \lim_{\varepsilon \rightarrow 0}
 \bar{\delta}^{x}_{\varepsilon^{-1}} \, \delta^{x}_{\varepsilon} u $$
Under closer scrutiny, it appears that the existence of the limit $\displaystyle 
Q^{x}$ (as a uniform limit, as well as having some other algebraic properties) is the basis 
which can be used for establishing sub-riemannian geometry.

\subsection*{Outline of the paper} After the introductory section \ref{slen} dedicated to 
basic notions concerning length in metric spaces, in section \ref{dilst} we 
describe the notion of a dilatation structure, introduced in \cite{buligadil1}. 
A dilatation structure on a metric space directly provides a notion of 
derivative, thus endowing the space with its own differential calculus. The
class of metric spaces admitting dilatation structures seems rather large,
containing riemannian, sub-riemannian as well as some ultrametric spaces, as
explained in \cite{buligadil2}, \cite{buligasr}, \cite{buligaultra}. The idea of
dilatation structures is that dilatations (or dilations, or homotheties, 
or even contractions as considered in the case of contractible groups) are
central objects for a differential calculus. The field $\delta$ 
of dilatations on a metric space $(X,d)$  obeys 5 axioms, 
see definition \ref{defweakstrong}, stating algebraic and analytical properties
of $\delta$, as well as the compatibility between $\delta$ and the distance $d$.

In section \ref{seclds} we propose an alternative notion, length dilatation
structures, which will be central in further considerations. In a length 
dilatation structure, definition \ref{deflds}, the accent is put on 
the length functional induced by the  distance $d$. We may imagine the field of dilatations 
$$\displaystyle (x, \varepsilon) \in X \times (0,1] \mapsto 
\delta^{x}_{\varepsilon} : U(x)\subset X \rightarrow X$$  as a field of 
microscopes with  magnification power $\varepsilon$, associating to any 
$x \in X$ a chart  $U(x)$ of a $\varepsilon$-neighbourhood of $x$, as measured
in the distance $d$. Imagine a curve in $X$ as a road and the various charts
provided by dilatations as roadmaps. In a length dilatation structure the
lengths of the images of the true road, as seen in different roadmaps, have to
agree. Also, these roadmaps have to be compatible in a clearly stated manner. 
Finally, the compatibility of the dilatation field with the length functional
induced by the distance $d$ is further stated as a Gamma-convergence condition
of induced length functionals, as $\varepsilon \rightarrow 0$.   

In  section \ref{sprop} is explained the structure of the tangent bundle which comes with 
a strong dilatation structure or  a length dilatation structure. The
characterization of the tangent bundle for length dilatation structures is new. 
A key notion which appears is the one of a conical group, studied in 
\cite{buligadil2}, which generalizes Carnot groups and contractible groups as 
well. 
 
In order to
facilitate the understanding of the abstract theory of tempered dilatation 
structures and coherent projections 
(sections \ref{stemp}, \ref{cohp} and \ref{sechormander}), we explain in section \ref{srsmooth} the case of dilatation structures
on sub-riemannian manifolds, following \cite{buligasr}. 

In section \ref{radon} we begin to study dilatation structures satisfying the
Radon-Nikodym property for metric spaces (or rectifiability
 property, or RNP), definition \ref{defrn}. This property says that 
 Lipschitz curves are derivable almost everywhere in the sense provided by 
 the dilatation structure. We give examples, then we easily obtain a description
 of the length functional as if we were in a kind of a generalized Finsler
 manifold, theorem \ref{fleng}. 
 
 Tempered dilatation structures, section \ref{stemp}, seem to be the habitat 
 where generalizations of results of Buttazzo, De Pascale and Fragala \cite{buttazzo1} and 
 Venturini \cite{venturini} naturally live. A dilatation structure is tempered, 
definition  \ref{dtempered}, if the charts provided by dilatations are
bi-lipschitz with the real distance, in a uniform manner with respect to 
the magnification $\varepsilon$ and the base point $x$. This is locally the
case for any $\displaystyle \mathcal{C}^{1}$ riemannian manifold, but it is not 
true for sub-riemannian manifolds, for example. From corollary \ref{cortemp} 
to theorem \ref{new3.1} we find out that a tempered dilatation structure with
RNP is also a length dilatation structure. 

In section \ref{cohp} coherent projections are introduced and studied. Coherent
projections are generalizations of distributions. With the help of a coherent 
projection $Q$ and a tempered dilatation structure $(X, \bar{d}, \bar{\delta})$ 
 we get a new field of dilatations $\delta$ and a new distance $d$, quite
 similar to a Carnot-Carath\'eodory distance. Notice however that in the case of
 sub-riemannian manifold we use as a tempered dilatation structure the one
 coming from a riemannian manifold, which according to our language has two 
 very special properties: it is locally linear (see the paper 
 \cite{buligadil2} for the affine geometry of a linear dilatation structure) 
 and it is commutative in the sense that the tangent spaces are commutative 
 conical groups, that is they are vector spaces. In the general formalism of 
 coherent projections and tempered dilatation structures nothing like this is
 used. 
 
 The main problem that we solve, section \ref{sechormander}, is if $(X,d, \delta)$ 
 is a length dilatation structure. This problem is solved for coherent projections which
 satisfy a generalized Chow condition. This condition  is inspired by the
 classical Chow condition, but for the reader which becomes familiar with
 dilatation structures is rather clear that Chow condition is only one among an
 infinity of other conditions with equivalent  effect. Indeed, even if we shall
 not touch this in the present paper, the Chow condition seems to be only a
 convenient way to indicate an algorithm for going from point A to point B, in
 terms of vector field brackets. We explained in \cite{buligadil1} that to
 dilatation structures in general is associated a formalism of binary decorated
 planar trees. At the level of this formalism the  algorithm from Chow condition, 
 as formulated in this paper, appears as working on a very particular class  
 of such binary trees. 
 
 In  subsection \ref{subscls} we finally get that coherent projections which
 satisfy condition (Cgen) and tempered dilatation structures which satisfy some 
 supplementary conditions (A) and (B) indeed induce length dilatation
 structures. At the classical level, this implies the new result that on regular
 sub-riemannian manifolds the rescaled (with the magnification factor
 $\varepsilon$) lengths Gamma-converge to the length in the metric tangent
 space, for any point. 
  
  The paper ends with the conclusion section \ref{seccon}, where we state 
  that Gromov'  viewpoint, that the CC distance is the only intrinsic object in sub-riemannian
  geometry, should be supplemented with Siebert' result, that a homogeneous Lie
  group is just a locally compact group endowed with a contracting and
  continuous  one parameter group of  automorphisms. This is  what we do 
  in this paper, by replacing the classical differential structures with 
  the more general dilatation structures. 


\section{Length in metric spaces}
\label{slen}

For a detailed introduction into the subject see for example \cite{amb}, chapter
1. 

\begin{definition}
The {\bf (upper) dilatation of a map} $f: X \rightarrow Y$ between metric spaces,  
in a point $u \in Y$ is 
$$ Lip(f)(u) = \limsup_{\varepsilon \rightarrow 0} \  
\sup  \left\{ 
\frac{d_{Y}(f(v), f(w))}{d_{X}(v,w)} \ : \ v \not = w \ , \ v,w \in B(u,\varepsilon)
 \right\}$$
\end{definition}
In the particular case of a   derivable function 
$f: \mathbb{R} \rightarrow \mathbb{R}^{n}$ the upper dilatation is  
$\displaystyle Lip(f)(t)  =  \|\dot{f}(t)\|$. 

 A function  $f:(X, d) \rightarrow (Y, d')$  is Lipschitz if there is a positive 
 constant $C$ such that for any $x,y \in X$ we have 
 $\displaystyle d'(f(x),f(y)) \leq C \, d(x,y)$. The number $Lip(f)$ is the smallest 
 such positive constant. Then  for any $x \in X$ we have the obvious relation  
$\displaystyle Lip(f)(x) \ \leq \ Lip(f)$.

A curve is a continuous function $c: [a,b] \rightarrow X$. The image of a curve is 
called path. Length measures paths. Therefore length does not depends on the 
reparameterization of the path and it is additive with respect to concatenation of paths.

\begin{definition}
In a metric space $(X,d)$ there are several ways to define  the length: 
\begin{enumerate}
\item[(a)] The {\bf length 
of a  curve with $L^{1}$ upper dilatation} $c: [a,b] \rightarrow X$ is 
$$L(f) = \int_{a}^{b} Lip(c)(t) \mbox{ d}t$$
\item[(b)] The {\bf variation of a curve}  $c: [a,b] \rightarrow X$ is the quantity 
$ Var(c)  =$  
$$  =  \sup \left\{ \sum_{i=0}^{n} d(c(t_{i}), 
c(t_{i+1})) \ \mbox{ : }  a = t_{0} <  t_{1} < ... < t_{n} < t_{n+1} =  b
\right\}$$ 
\item[(c)] The {\bf length of the path} $A = c([a,b])$ is  the
one-dimensional Hausdorff measure of the path.:  
$$l(A) \ = \ \lim_{\delta \rightarrow 0}  
 \inf \left\{ \sum_{i \in I} diam \ E_{i}  \mbox{ : } diam \ E_{i} 
< \delta \ , \ \ A \subset \bigcup_{i \in I} E_{i} \right\} $$
\end{enumerate}
\label{deflenght}
\end{definition}

The definitions are not equivalent. For Lipschitz curves the first  two 
definitions agree. For simple Lipschitz  curves all definitions agree.

\begin{theorem}
For each Lipschitz curve $c: [a,b] \rightarrow X$, we have 
$$\displaystyle L(c) \ = \ Var(c) \ \geq \ \mathcal{H}^{1}(c([a,b]))$$

 If $c$ is moreover injective then  $\displaystyle \mathcal{H}^{1}(c([a,b])) \  = 
 \ Var(f)$. 
\label{t411amb}
\end{theorem}

An important tool used in the proof  of the previous theorem is 
the geometrically obvious, but not straightforward to prove in this generality, 
reparametrisation Theorem. 

\begin{theorem}
Any Lipschitz curve  $c$ admits a reparametrisation $c'$ such that $Lip(c')(t) = 1$ for almost any $t \in [a,b]$. 
\label{tp}
\end{theorem}

\begin{definition}
We shall denote by $l_{d}$ the {\bf length functional induced by the distance} $d$, 
defined only on the family of Lipschitz curves. 
  If the metric space $(X,d)$ is connected by Lipschitz curves, then the length 
induces a new distance $d_{l}$, given by: 
$$d_{l}(x,y) \  = \ \inf \ \left\{ l_{d}(c([a,b])) \mbox{ : } 
c: [a,b] \rightarrow X \ \mbox{ Lipschitz } , \right.$$
$$\left. \ c(a)=x \ , \ c(b) = y \right\}$$

A {\bf length metric space} is a metric space $(X,d)$, connected by Lipschitz curves, 
 such that $d  = d_{l}$. 
\label{dpath}
\end{definition}

From theorem \ref{t411amb} we deduce that Lipschitz curves in complete 
length metric spaces are absolutely continuous. 
Indeed, here is the definition of an absolutely continuous curve 
(definition 1.1.1, chapter 1,   \cite{amb}). 
 
 \begin{definition}
 Let $(X,d)$ be a complete metric space. A curve $c:(a,b)\rightarrow X$ is {\bf absolutely 
 continuous} if there exists $m\in L^{1}((a,b))$ such that for any $a<s\leq t<b$ we have 
 $$d(c(s),c(t)) \leq \int_{s}^{t} m(r) \mbox{ d}r   .$$
 Such a function $m$ is called a {\bf upper gradient} of the curve $c$. 
 \label{defac}
 \end{definition}
 
 According to theorem \ref{t411amb}, for a Lipschitz curve $c:[a,b]\rightarrow X$ in a 
 complete length metric space such a function 
 $m\in L^{1}((a,b))$  is the upper dilatation  $Lip(c)$. 
More can be said about the expression of the upper dilatation. We need first to introduce the notion of 
metric derivative of a Lipschitz curve. 

\begin{definition}
A curve $c:(a,b)\rightarrow X$ is {\bf metrically derivable} in $t\in(a,b)$ if the limit 
$$md(c)(t) = \lim_{s\rightarrow t} \frac{d(c(s),c(t))}{\mid s-t \mid}$$
exists and it is finite. In this case $md(c)(t)$ is called the {\bf metric
derivative} of $c$ in $t$. 
\label{defmd}
\end{definition}

For the proof of the following theorem see \cite{amb}, theorem 1.1.2, chapter 1. 

\begin{theorem}
Let $(X,d)$ be a complete metric space and $c:(a,b)\rightarrow X$ be an absolutely continuous curve. 
Then $c$ is metrically  derivable for $\mathcal{L}^{1}$-a.e. $t\in(a,b)$. Moreover the function $md(c)$ belongs to $L^{1}((a,b))$ and it is minimal in the following sense: $md(c)(t)\leq m(t)$  for  $\mathcal{L}^{1}$-a.e. $t\in(a,b)$, for each upper gradient $m$ of the curve $c$. 
\label{tupper}
\end{theorem}

\section{Dilatation structures}
\label{dilst}

We shall use here a slightly particular version of dilatation structures. 
For the general definition of a dilatation structure see \cite{buligadil1} (the general definition 
applies for dilatation structures over ultrametric spaces as well).

\begin{definition}
Let $(X,d)$ be a complete metric space such that for any $x  \in X$ the 
closed ball $\bar{B}(x,3)$ is compact. A {\bf dilatation structure} $(X,d, \delta)$ 
over $(X,d)$ is the assignment to any $x \in X$  and $\varepsilon \in (0,+\infty)$ 
of a invertible homeomorphism, defined as: if 
$\displaystyle   \varepsilon \in (0, 1]$ then  $\displaystyle 
 \delta^{x}_{\varepsilon} : U(x)
\rightarrow V_{\varepsilon}(x)$, else 
$\displaystyle  \delta^{x}_{\varepsilon} : 
W_{\varepsilon}(x) \rightarrow U(x)$,  such that the following axioms are satisfied: 
\begin{enumerate}
\item[{\bf A0.}]  For any $x \in X$ the sets $ \displaystyle U(x), V_{\varepsilon}(x), 
W_{\varepsilon}(x)$ are open neighbourhoods of $x$.  There are  numbers  $1<A<B$ such that for any $x \in X$  and any 
$\varepsilon \in (0,1)$ we have 
  the following string of inclusions:
$$ B_{d}(x, \varepsilon) \subset \delta^{x}_{\varepsilon}  B_{d}(x, A) 
\subset V_{\varepsilon}(x) \subset 
W_{\varepsilon^{-1}}(x) \subset \delta_{\varepsilon}^{x}  B_{d}(x, B) $$
Moreover for any compact set $K \subset X$ there are $R=R(K) > 0$ and 
$\displaystyle \varepsilon_{0}= \varepsilon(K) \in (0,1)$  such that  
for all $\displaystyle u,v \in \bar{B}_{d}(x,R)$ and all 
$\displaystyle \varepsilon  \in (0,\varepsilon_{0})$,  we have 
$$\delta_{\varepsilon}^{x} v \in W_{\varepsilon^{-1}}( \delta^{x}_{\varepsilon}u) \ .$$

\item[{\bf A1.}]  We  have 
$\displaystyle  \delta^{x}_{\varepsilon} x = x $ for any point $x$. 
We also have $\displaystyle \delta^{x}_{1} = id$ for any $x \in X$. 
Let us define the topological space
$$ dom \, \delta = \left\{ (\varepsilon, x, y) \in (0,+\infty) \times X 
\times X \mbox{ :  if } \varepsilon \leq 1 \mbox{ then } y 
\in U(x) \,
\, , 
\right.$$ 
$$\left. \mbox{  else } y \in W_{\varepsilon}(x) \right\} $$ 
with the topology inherited from $(0,+\infty) \times X \times X$ endowed with
 the product topology. Consider also 
$\displaystyle Cl(dom \, \delta)$, 
the closure of 
$dom \, \delta$ in $\displaystyle [0,+\infty) \times X \times X$. The function $\displaystyle \delta : dom \, \delta 
\rightarrow  X$ defined by $\displaystyle \delta (\varepsilon,  x, y)  = 
\delta^{x}_{\varepsilon} y$ is continuous. Moreover, it can be continuously 
extended to the set $\displaystyle Cl(dom \, \delta)$ and we have 
$$\lim_{\varepsilon\rightarrow 0} \delta_{\varepsilon}^{x} y \, = \, x  $$

\item[{\bf A2.}] For any  $x, \in X$, $\displaystyle \varepsilon, \mu \in (0,+\infty)$
 and $\displaystyle u \in U(x)$   we have the equality: 
$$ \delta_{\varepsilon}^{x} \delta_{\mu}^{x} u  = \delta_{\varepsilon \mu}^{x} u $$ 
whenever one of the sides are well defined.

\item[{\bf A3.}]  For any $x$ there is a distance  function $\displaystyle (u,v) \mapsto d^{x}(u,v)$, defined for any $u,v$ in the closed ball (in distance d) $\displaystyle 
\bar{B}(x,A)$, such that 
$$\lim_{\varepsilon \rightarrow 0} \quad \sup  \left\{  \mid 
\frac{1}{\varepsilon} d(\delta^{x}_{\varepsilon} u, \delta^{x}_{\varepsilon} v) \ - \ d^{x}(u,v) \mid \mbox{ :  } u,v \in \bar{B}_{d}(x,A)\right\} \ =  \ 0$$
uniformly with respect to $x$ in compact set. 

\end{enumerate}

The {\bf dilatation structure is strong} if it satisfies the following supplementary condition: 

\begin{enumerate}
\item[{\bf A4.}] Let us define 
$\displaystyle \Delta^{x}_{\varepsilon}(u,v) =
\delta_{\varepsilon^{-1}}^{\delta^{x}_{\varepsilon} u} \delta^{x}_{\varepsilon} v$. 
Then we have the limit 
$$\lim_{\varepsilon \rightarrow 0}  \Delta^{x}_{\varepsilon}(u,v) =  \Delta^{x}(u, v)  $$
uniformly with respect to $x, u, v$ in compact set. 
\end{enumerate}
\label{defweakstrong}
\end{definition}
 
We shall use many times from now the words "sufficiently close". This deserves 
a definition. 

\begin{definition}
Let  $(X,d, \delta)$ be a strong dilatation structure. We say that a property 
$$\displaystyle \mathcal{P}(x_{1},x_{2},
x_{3}, ...)$$ is true  for $\displaystyle x_{1}, x_{2}, x_{3},
...$ {\bf sufficiently close} if for any compact, non empty set $K \subset X$, there
is a positive constant $C(K)> 0$ such that $\displaystyle \mathcal{P}(x_{1},x_{2},
x_{3}, ...)$ is true for any $\displaystyle x_{1},x_{2},
x_{3}, ... \in K$ with $\displaystyle d(x_{i}, x_{j}) \leq C(K)$.
\end{definition}
 
\section{Length dilatation structures}
\label{seclds}

Consider $(X,d)$ a complete, locally compact metric space, and a triple 
 $(X,d,\delta)$  which satisfies  A0, A1, A2. Denote by $\displaystyle Lip([0,1],X,d)$ the space of 
$d$-Lipschitz curves $c:[0,1] \rightarrow X$. Let also $\displaystyle 
l_{d}$ denote the length functional associated to the distance $d$.

\subsection*{Gamma-convergence of length functionals}

\begin{definition}
For any $\varepsilon \in (0,1)$ we define the {\bf length functional }
$$l_{\varepsilon}: \mathcal{L}_{\varepsilon}(X,d,  \delta) \rightarrow
[0,+\infty] \quad , \quad l_{\varepsilon}(x,c) \ = \ l^{x}_{\varepsilon}(c) \ = \ \frac{1}{\varepsilon}
\, 
l_{d}(\delta^{x}_{\varepsilon} c)  $$
The domain of definition of the functional 
$\displaystyle l_{\varepsilon}$ is the space: 
 $$\mathcal{L}_{\varepsilon}(X,d,  \delta) \ = \ \left\{(x ,c) \in X  
 \times \mathcal{C}([0,1],X) 
\mbox{ : } c: [0,1] \in U(x) \, \, , \,  \right.$$ 
$$\left. \delta^{x}_{\varepsilon}c \mbox{ is } 
d-Lip \mbox{ and } Lip(\delta^{x}_{\varepsilon}c)\,  \leq \, 2  \, 
l_{d}(\delta^{x}_{\varepsilon}c) 
\right\} $$
\label{thespaceleps}
\end{definition}

The last condition from the definition of $\displaystyle 
\mathcal{L}_{\varepsilon}(X,d,  \delta)$ is a selection of parameterization 
of the path $c([0,1])$. Indeed, by the reparameterization theorem, if 
$\displaystyle \delta^{x}_{\varepsilon} c :[0,1] \rightarrow (X,d)$ is a 
$d$-Lipschitz curve of length $\displaystyle L = l_{d}(\delta^{x}_{\varepsilon}c)$ 
then $\displaystyle \delta^{x}_{\varepsilon}c([0,1])$ can be reparameterized by length, that is there exists a 
increasing  function $\phi:[0,L] \rightarrow [0,1]$ such that $\displaystyle 
c'= \delta^{x}_{\varepsilon} c\circ \phi$ is 
a $d$-Lipschitz curve with $Lip(c') \leq 1$. But we can use a second affine
reparameterization which sends $[0,L]$ back to $[0,1]$ and we get a Lipschitz
curve $c"$ with $c"([0,1]) = c'([0,1])$ and $\displaystyle Lip(c") \leq 2 
l_{d}(c)$.

We shall use the following definition of Gamma-convergence (see the book
 \cite{dalmaso} for the notion of Gamma-convergence). Notice the use of
convergence of sequences only in the second part of the definition.

\begin{definition}
Let $Z$ be a metric space with distance function $D$ and $\displaystyle \left(
l_{\varepsilon}\right)_{\varepsilon > 0}$ be a family of functionals $\displaystyle 
l_{\varepsilon}: Z_{\varepsilon} \subset Z \rightarrow [0,+\infty]$. Then 
$\displaystyle l_{\varepsilon}$ {\bf Gamma-converges} to the functional 
$\displaystyle l: Z_{0} \subset Z \rightarrow [0,+\infty]$ if: 
\begin{enumerate}
\item[(a)] ({\bf liminf inequality}) for any function $\displaystyle \varepsilon \in
(0,\infty)  \mapsto 
x_{\varepsilon} \in Z_{\varepsilon}$ such that $\displaystyle \lim_{\varepsilon
\rightarrow 0} x_{\varepsilon} \, = \, x_{0} \in Z_{0}$ we have 
$$l(x_{0}) \, \leq \, \liminf_{\varepsilon \rightarrow 0}
l_{\varepsilon}(x_{\varepsilon})$$
\item[(b)] ({\bf existence of a recovery sequence}) For any $\displaystyle x_{0} \in Z_{0}$ 
and for any sequence $\displaystyle \left( \varepsilon_{n} \right)_{n \in \mathbb{N}}$
such that $\displaystyle \lim_{n \rightarrow \infty} \varepsilon_{n} \, = \, 0$ there
is a sequence $\displaystyle \left( x_{n} \right)_{n \in \mathbb{N}}$ with
$\displaystyle x_{n} \in Z_{\varepsilon_{n}}$ for any $n \in \mathbb{N}$, such that 
$$l(x_{0}) \, = \, \lim_{n \rightarrow \infty}
l_{\varepsilon_{n}}(x_{n})$$
\end{enumerate}
\end{definition}

We shall take as the metric space $Z$ the space
 $X \times \mathcal{C}([0,1],X)$ with the distance 
 $$ D((x,c), (x', c')) \ = \ \max\left\{ d(x, x') \, , \, \sup \left\{ d(c(t),
 c'(t)) \mbox{ : } t \in [0,1] \right\} \right\} $$

Let  $\displaystyle \mathcal{L}(X,d,\delta)$be the class of all 
$\displaystyle (x ,c) \in X  \times \mathcal{C}([0,1],X)$ which appear as  limits  
$\displaystyle (x_{n}, c_{n}) \rightarrow (x,c)$, with $\displaystyle 
(x_{n}, c_{n}) \in \mathcal{L}_{\varepsilon_{n}}(X,d,  \delta)$, 
the family $\displaystyle (c_{n})_{n}$ is $d$-equicontinuous and $\displaystyle 
\varepsilon_{n} \rightarrow 0$ as $n \rightarrow \infty$.

\begin{definition}
A triple $(X,d,\delta)$ is a {\bf length dilatation structure} if $(X,d)$ is a
complete, locally compact metric space such that A0, A1, 
A2,  are satisfied, together with  the following axioms: 
\begin{enumerate}
\item[\bf{A3L.}] there is a functional $\displaystyle l : \mathcal{L}(X,d,  \delta) \rightarrow
[0,+\infty]$ such that for any $\varepsilon_{n} \rightarrow 0$ as $n \rightarrow
\infty$ the sequence of functionals 
$\displaystyle l_{\varepsilon_{n}}$ Gamma-converges to the functional $l$. 
\item[\bf{A4+}] Let us define 
$\displaystyle \Delta^{x}_{\varepsilon}(u,v) =
\delta_{\varepsilon^{-1}}^{\delta^{x}_{\varepsilon} u} \delta^{x}_{\varepsilon}
v$ and $\displaystyle \Sigma^{x}_{\varepsilon}(u,v) = 
\delta_{\varepsilon^{-1}}^{x} \delta_{\varepsilon}^{\delta^{x}_{\varepsilon} u}
v$. 
Then we have the limits 
$$\lim_{\varepsilon \rightarrow 0}  \Delta^{x}_{\varepsilon}(u,v) =  \Delta^{x}(u, v)  $$
$$\lim_{\varepsilon \rightarrow 0}  \Sigma^{x}_{\varepsilon}(u,v) =  \Sigma^{x}(u, v)  $$
uniformly with respect to $x, u, v$ in compact set. 
\end{enumerate} 
\label{deflds}
\end{definition}

\begin{remark}
For strong dilatation structures the axioms A0 - A4 imply A4+, cf. corollary 9  \cite{buligadil1}. The
transformations $\displaystyle \Sigma^{x}_{\varepsilon}(u, \cdot)$ have the
interpretation of approximate left translations in the tangent space of $(X,d)$
at $x$.  
\end{remark}

For any $\varepsilon \in (0,1)$ and any $x \in X$ the length functional 
$\displaystyle l^{x}_{\varepsilon}$ induces a distance on $U(x)$: 
$$ \mathring{d}^{x}_{\varepsilon}(u,v) \ = \ \inf\left\{ l^{x}_{\varepsilon}(c)
\mbox{ : } (x,c) \in \mathcal{L}_{\varepsilon}(X,d,  \delta) \, , \, c(0) = u \,
, \, c(1) = v \right\} $$
In the same way the length functional $l$ from A3L induces a distance
$\displaystyle \mathring{d}^{x}$ on $U(x)$.

Gamma-convergence implies that 
\begin{equation}
\mathring{d}^{x}(u,v) \, \geq \, \limsup_{\varepsilon \rightarrow 0} \mathring{d}^{x}_{\varepsilon}(u,v)
\label{dsup}
\end{equation}

\begin{remark}
Without supplementary hypotheses we cannot prove A3 from A3L, that is 
in principle length dilatation structures are not strong dilatation structures. 
\label{rkused}
\end{remark}

\section{Properties of (length) dilatation structures}
\label{sprop}

For a dilatation structure the metric tangent spaces   have a group structure which is 
compatible with dilatations.

 We shall work further with local groups. Such objects are spaces endowed with 
a locally defined  operation,  satisfying the conditions of a uniform group. 
 See section 3.3 \cite{buligadil1} for details about the 
definition of local groups.

\subsection{Normed conical groups}

These have been introduced in section 8.2 \cite{buligadil1} and studied further in section 4 
\cite{buligadil2}. In the following general definition appear a topological commutative group 
$\Gamma$ endowed with a continuous morphism $\nu: \Gamma \rightarrow (0, +\infty)$ from $\Gamma$ 
to the group $(0, +\infty)$ with multiplication.  The morphism $\nu$ induces an invariant topological filter on $\Gamma$ (other names for such an invariant filter are "absolute" or "end").  The convergence 
of a variable $\varepsilon \in \Gamma$ to this filter is denoted by $\varepsilon \rightarrow 0$ and 
it means simply $\nu(\varepsilon) \rightarrow 0$ in $\mathbb{R}$. 

Particular, interesting examples of pairs $(\Gamma, \nu)$ are: $(0, +\infty)$ with identity, which is 
the case interesting for this paper, $\displaystyle \mathbb{C}^{*}$ with the modulus of complex numbers, 
or $\mathbb{N}$ (with addition) with the exponential, which is relevant for the case of normed 
contractible groups, section 4.3 \cite{buligadil2}. 

\begin{definition}
A {\bf normed group with dilatations} $(G, \delta, \| \cdot \|)$ is a local  group 
$G$  with  a local action of $\Gamma$ (denoted by $\delta$), on $G$ such that
\begin{enumerate}
\item[H0.] the limit  $\displaystyle \lim_{\varepsilon \rightarrow 0}
\delta_{\varepsilon} x  =  e$ exists and is uniform with respect to $x$ in a compact neighbourhood of the identity $e$.
\item[H1.] the limit
$$\beta(x,y)  =  \lim_{\varepsilon \rightarrow 0} \delta_{\varepsilon}^{-1}
\left((\delta_{\varepsilon}x) (\delta_{\varepsilon}y ) \right)$$
is well defined in a compact neighbourhood of $e$ and the limit is uniform.
\item[H2.] the following relation holds
$$ \lim_{\varepsilon \rightarrow 0} \delta_{\varepsilon}^{-1}
\left( ( \delta_{\varepsilon}x)^{-1}\right)  =  x^{-1} $$
where the limit from the left hand side exists in a neighbourhood of 
$e$ and is uniform with respect to $x$.
\end{enumerate}

Moreover the group is endowed with a continuous norm
function $\displaystyle \|\cdot \| : G \rightarrow \mathbb{R}$ which satisfies
(locally, in a neighbourhood  of the neutral element $e$) the properties:
 \begin{enumerate}
 \item[(a)] for any $x$ we have $\| x\| \geq 0$; if $\| x\| = 0$ then $x=e$,
 \item[(b)] for any $x,y$ we have $\|xy\| \leq \|x\| + \|y\|$,
 \item[(c)] for any $x$ we have $\displaystyle \| x^{-1}\| = \|x\|$,
 \item[(d)] the limit
$\displaystyle \lim_{\varepsilon \rightarrow 0} \frac{1}{\nu(\varepsilon)} \| \delta_{\varepsilon} x \| = \| x\|^{N}$
 exists, is uniform with respect to $x$ in compact set,
 \item[(e)] if $\displaystyle \| x\|^{N} = 0$ then $x=e$.
  \end{enumerate}
  \label{dnco}
  \end{definition}

In a normed group with dilatations we have a natural left invariant distance 
given by
\begin{equation}
d(x,y) = \| x^{-1}y\| 
\label{dnormed}
\end{equation}
Any locally compact  normed group with dilatations has an associated dilatation structure on it.  In a group with dilatations $(G, \delta)$  we define dilatations based in any point $x \in G$ by
 \begin{equation}
 \delta^{x}_{\varepsilon} u = x \delta_{\varepsilon} ( x^{-1}u)  .
 \label{dilat}
 \end{equation}

The following result is theorem 15 \cite{buligadil1}.

\begin{theorem}
Let $(G, \delta, \| \cdot \|)$ be  a locally compact  normed local group with dilatations. Then $(G, d, \delta)$ is
a strong dilatation structure, where $\delta$ are the dilatations defined by (\ref{dilat}) and the distance $d$ is induced by the norm as in (\ref{dnormed}).
\label{tgrd}
\end{theorem}

 \begin{remark}
 In \cite{buligadil1} dilatation structures are also defined using a pair 
 $(\Gamma, \nu)$ and the convention of writing $\varepsilon \rightarrow 0$, 
 $\varepsilon \in \Gamma$, instead of $\nu(\varepsilon) \rightarrow 0$. 
 In this paper $\Gamma$ will always be taken as equal to $(0,+\infty)$, but 
 definitions and results described further can be expressed for a general pair 
 $(\Gamma, \nu)$. 
 \end{remark}

\begin{definition}
A {\bf normed conical group} $N$ is a normed  group with dilatations  such that
for any $\varepsilon \in \Gamma$  the dilatation
 $\delta_{\varepsilon}$ is a group morphism  and such that for any $\varepsilon >0$
  $\displaystyle  \| \delta_{\varepsilon} x \| = \nu(\varepsilon) \| x \|$.
\end{definition}

A conical group is the infinitesimal version of a group with
dilatations (\cite{buligadil1} proposition 2).

\begin{proposition}
Under the hypotheses H0, H1, H2 $\displaystyle (G,\beta, \delta, \| \cdot \|^{N})$ is a locally compact, 
local normed conical group, with operation $\beta$,  dilatations $\delta$ and homogeneous norm $\displaystyle \| \cdot \|^{N}$.
\label{here3.4}
\end{proposition}

\subsection{Tangent bundle of a dilatation structure}
\label{induced}

The following two theorems describe the most important metric and algebraic 
properties of a dilatation structure. As presented 
here these are condensed statements, available in full length as theorems 7, 8,
10 in \cite{buligadil1}. The first theorem does not need a proof (see theorem 7 
\cite{buligadil1}).

\begin{theorem}
Let $(X,d,\delta)$ be a  strong dilatation structure. Then the metric space $(X,d)$ 
admits a metric tangent space at $x$, for any point $x\in X$. 
More precisely we have  the following limit: 
$$\lim_{\varepsilon \rightarrow 0} \ \frac{1}{\varepsilon} \sup \left\{  \mid d(u,v) - d^{x}(u,v) \mid \mbox{ : } d(x,u) \leq \varepsilon \ , \ d(x,v) \leq \varepsilon \right\} \ = \ 0 \ .$$
\label{thcone}
\end{theorem}

 Length dilatation structures were introduced in this
paper. Straightforward modifications in the proof of the before mentioned
theorems allow us to extend some results to  length dilatation structures.

\begin{theorem}
Let $(X,d,\delta)$ be a strong dilatation structure or a length dilatation
structure. Then: 
 \begin{enumerate}
 \item[(a)]  $\displaystyle \Sigma^{x}$ is a local group operation on $U(x)$,
 with $x$ as neutral element and $\displaystyle \, inv^{x}$ as the inverse element
 function; 
  \item[(b)] for strong dilatation structures the distance 
  $\displaystyle d^{x}$ is left invariant with respect 
  to the group operation from point (a); for length dilatation structures the
  length functional $\displaystyle l^{x} = l(x, \cdot)$ is invariant with
  respect to left translations $\displaystyle \Sigma^{x}(y, \cdot)$, $y \in
  U(x)$;  
 \item[(c)] For any $\varepsilon \in (0,1]$ the 
 dilatation $\displaystyle \delta^{x}_{\varepsilon}$ is an automorphism with 
 respect to the group operation from point (a); 
 \item[(d)] for strong dilatation structures the distance $d^{x}$ has the cone property with
respect to dilatations: for any $u,v \in X$ such that $\displaystyle d(x,u)\leq 1$ and 
$\displaystyle d(x,v) \leq 1$  and all $\mu \in (0,A)$ we have: 
$$d^{x}(u,v) \ = \ \frac{1}{\mu} d^{x}(\delta_{\mu}^{x} u , \delta^{x}_{\mu} v) 
 $$ 
 For length dilatation structures we have for any $\mu \in (0,1]$ the equality
 $$l(x, \delta^{x}_{\mu} c) \, =  \, \mu \, l (x, c)$$
 \end{enumerate}
\label{tgene}
\end{theorem}

\paragraph{Proof.}
We shall only prove the statements concerning length dilatation structures. For
(a) and (c) notice that the axiom A4+ is all that we need in order to
transform  the proof of theorem 10 \cite{buligadil1} into a proof of this point.
Indeed, for this we need the existence of the limits from A4+ and the algebraic
relations from theorem 11 \cite{buligadil1} which are true only from A0, A1, A2.

For (b) remark that if $\displaystyle (\delta^{x}_{\varepsilon}y ,c) \in
\mathcal{L}_{\varepsilon}(X,d,\delta)$ then $\displaystyle 
(x, \Sigma^{x}_{\varepsilon}(y, \cdot) c) \in
\mathcal{L}_{\varepsilon}(X,d,\delta)$ and moreover 
$$l_{\varepsilon}(\delta^{x}_{\varepsilon}y ,c) \, = \, l_{\varepsilon} 
(x, \Sigma^{x}_{\varepsilon}(y, \cdot) c)$$
Indeed, this is true because of the equality: 
$$\delta^{\delta^{x}_{\varepsilon} y}  c \, = \, \delta^{x}_{\varepsilon} 
\Sigma^{x}_{\varepsilon}(y, \cdot) c$$
By passing to the limit with $\varepsilon \rightarrow 0$ and using A3L and A4+ we
get 
$$l(x, c) \, = \, l(x, \Sigma^{x}(y, \cdot) c)$$

For (d) remark that for any $\varepsilon, \mu > 0$ (and sufficiently small) 
$\displaystyle (x,c) \in \mathcal{L}_{\varepsilon \mu}(X,d,\delta)$ is equivalent with 
$\displaystyle (x, \delta^{x}_{\mu} c) \in
\mathcal{L}_{\varepsilon}(X,d,\delta)$ 
and moreover: 
$$l_{\varepsilon}(x, \delta^{x}_{\mu} c) \, = \, \frac{1}{\varepsilon} \, 
l_{d}(\delta^{x}_{\varepsilon \mu} c) \, = \, \mu \, l_{\varepsilon \mu} (x, c)$$
We pass to the limit with $\varepsilon \rightarrow 0$ and we get the desired
equality.
\hfill $\square$

\begin{definition}
The conical group $\displaystyle (U(x), \Sigma^{x}, \delta^{x})$ can be seen as the tangent space 
of $(X,d, \delta)$ at $x$. We shall  denote it by  
$\displaystyle T_{x} (X, d, \delta) =  (U(x), \Sigma^{x}, \delta^{x})$, or by $\displaystyle T_{x} X$ if 
$(d,\delta)$ are clear from the context. 
\end{definition}

The following proposition is  corollary 
6.3 from \cite{buligadil2}, which gives a more precise
description of the conical group 
$\displaystyle (U(x), \Sigma^{x}, \delta^{x})$.

\begin{proposition}
Let $(X,d,\delta)$ be a strong dilatation structure. 
Then for any $x \in X$ the local group
 $\displaystyle (U(x), \Sigma^{x})$ is locally a simply connected Lie group
 whose Lie algebra admits a positive graduation (a homogeneous group), given by
 the eigenspaces of $\displaystyle \delta^{x}_{\varepsilon}$ for an 
 arbitrary $\varepsilon \in (0,1)$.
\label{cor63}
\end{proposition}

\begin{remark} 
S. Vodopyanov (private communication) made the observation that in the proof 
of  corollary 6.3  \cite{buligadil2} it is used  Siebert' proposition 5.4 \cite{siebert}, which 
is true for conical groups (in our language), while I am using it for {\em local}Ê   conical groups. 
This is true and constitutes a gap in the proof of the corollary 6.3. Fortunately the recent 
paper \cite{recent} provides the needed result for local groups. Indeed, theorem 1.1 \cite{recent} 
states that a locally compact, locally connected, contractible (with  Siebert' wording) group is 
locally isomorphic to a contractive Lie group.   
\end{remark} 

\subsection{Differentiability with respect to dilatation structures}

For any strong dilatation structure or length dilatation structure there is an associated  notion  of differentiability (section 7.2 \cite{buligadil1}). 
First we need the definition of a morphism of conical groups. 

\begin{definition}
 Let $(N,\delta)$ and $(M,\bar{\delta})$ be two  conical groups. A function $f:N\rightarrow M$ is a conical group morphism if $f$ is a group morphism and for any $\varepsilon>0$ and $u\in N$ we have 
 $\displaystyle f(\delta_{\varepsilon} u) = \bar{\delta}_{\varepsilon} f(u)$. 
\label{defmorph}
\end{definition}

The definition of the derivative, or differential,  with respect to dilatations structures follows. In the case of a pair of Carnot groups this is just the definition of the Pansu derivative 
introduced in \cite{pansu}.

 \begin{definition}
 Let $(X, d, \delta)$ and $(Y, \overline{d}, \overline{\delta})$ be two 
 strong dilatation structures  or length and $f:X \rightarrow Y$ be a continuous function. The function $f$ is differentiable in $x$ if there exists a 
 conical group morphism  $\displaystyle D \, f(x):T_{x}X\rightarrow T_{f(x)}Y$, defined on a neighbourhood of $x$ with values in  a neighbourhood  of $f(x)$ such that 
\begin{equation}
\lim_{\varepsilon \rightarrow 0} \sup \left\{  \frac{1}{\varepsilon} \overline{d} \left( f\left( \delta^{x}_{\varepsilon} u\right) ,  \overline{\delta}^{f(x)}_{\varepsilon} D \, f(x)  (u) \right) \mbox{ : } d(x,u) \leq \varepsilon \right\}Ê  = 0 , 
\label{edefdif}
\end{equation}
The morphism $\displaystyle D \, f(x) $ is called the derivative, or differential,  of $f$ at $x$.

\label{defdiffer}
\end{definition}

The definition also makes sense if the function $f$ is defined on a open subset of $(X,d)$. 

\section{Dilatation structures on sub-riemannian manifolds}
\label{srsmooth}

In \cite{buligasr}  we proved that we can associate 
  dilatation structures to regular
sub-Riemannian manifolds. This result, explained further, is the source of inspiration
of the notion of a coherent projection (section \ref{cohp}).

Let $M$ be a connected $n$ dimensional real manifold. A distribution
is a smooth  subbundle $D$ of $M$. To any point
$x \in M$ there is associated the vector space $D_{x} \subset T_{x}M$.
The dimension of the distribution $D$ at point $x \in M$ is denoted by
$$m(x) = \, dim \, D_{x}$$
The distribution is smooth, therefore the function $x \in M \mapsto m(x)$ is
locally constant. We suppose further that the dimension of the distribution is
globally constant and we  denote it by $m$ (thus $m = m(x)$ for any
$x \in M$). Clearly $m \leq n$; we are interested in the case $m < n$.

A horizontal curve $c:[a,b] \rightarrow M$ is a curve which is almost everywhere derivable and for
almost any $t \in [a,b]$ we have $\displaystyle \dot{c}(t) \in D_{c(t)}$.
The class of horizontal curves will be denoted by $Hor(M,D)$.

Further we shall use the following notion of non-integrability of the
distribution $D$.

\begin{definition}
The distribution $D$ is {\bf completely non-integrable} if $M$ is locally connected
by horizontal curves curves  $c \in Hor(M,D)$.
\label{defcnint}
\end{definition}

The Chow condition (C) \cite{chow} gives a sufficient condition for the 
distribution $D$ to be completely non-integrable.

\begin{theorem} (Chow) Let $D$ be a distribution of dimension $m$  in the manifold
$M$. Suppose there is a positive integer number $k$ (called the rank of the
distribution $D$) such that for any $x \in X$ there is a topological  open ball
$U(x) \subset M$ with $x \in U(x)$ such that there are smooth vector fields
$\displaystyle X_{1}, ..., X_{m}$ in $U(x)$ with the property:

(C) the vector fields $\displaystyle X_{1}, ..., X_{m}$ span $\displaystyle
D_{y}$ and these vector fields together with  their iterated
brackets of order at most $k$ span the tangent space $\displaystyle T_{y}M$
at every point $y \in U(x)$.

Then the distribution $D$ is completely non-integrable in the sense of
definition \ref{defcnint}.
\label{tchow}
\end{theorem}

\begin{definition}
A {\bf sub-riemannian  manifold} or SR manifold is a triple $(M,D, g)$, where $M$ is a
connected manifold, $D$ is a completely non-integrable distribution on $M$, 
and $g$ is a metric (Euclidean inner-product) on the distribution 
(or horizontal bundle)  $D$.
\label{defsr}
\end{definition}

Given a distribution $D$ which satisfies the hypothesis of Chow theorem
\ref{tchow}, let us consider a point   $x \in M$,  its neighbourhood $U(x)$,
and the vector fields $\displaystyle X_{1}, ..., X_{m}$ satisfying the
condition (C).

One can define on $U(x)$ a filtration of bundles
as follows. Define first the class of horizontal vector fields on $U$:
$$\mathcal{X}^{1}(U(x),D) \ = \ \left\{ X \in \mathcal{X}^{\infty}(U) \mbox{ : }
\forall y \in U(x) \ , \ X(y) \in D_{y} \right\}$$
Next, define inductively for all positive integers $j$:
$$ \mathcal{X}^{j+1} (U(x),D) \ = \ \mathcal{X}^{j}(U(x),D) \, + \,
 [ \mathcal{X}^{1}(U(x),D),
\mathcal{X}^{j}(U(x),D)]$$
Here $[ \cdot , \cdot ]$ denotes the bracket of vector fields.
We obtain therefore a filtration $\displaystyle \mathcal{X}^{j}(U(x),D) \subset
\mathcal{X}^{j+1} (U(x),D)$.
Evaluate now this filtration at $y \in U(x)$:
$$V^{j}(y,U(x),D) \ = \ \left\{ X(y) \mbox{ : } X \in \mathcal{X}^{j}(U(x),D)\right\}$$
According to Chow theorem there is a positive integer $k$ such that
 for all $y \in U(x)$ we have
$$D_{y} =  V^{1}(y,U(x),D) \subset V^{2}(y,U(x),D) \subset ... \subset
V^{k}(y, U(x),D) = T_{y}M $$
Consequently, to the sub-riemannian manifold is associated the string of
numbers:
$$ \nu_{1}(y) = \dim V^{1}(y, U(x),D) < \nu_{2}(y) = \dim V^{2}(y, U(x),D)
< ... < n = \dim M$$
Generally $k$, $\displaystyle \nu_{j}(y)$
may vary from a point to another.

The number $k$ is called the step of the distribution at
$y $.

\begin{definition}
The distribution $D$ is regular if  $\displaystyle \nu_{j}(y)$ are constant on
the manifold $M$.
The sub-riemannian manifold $(M,D,g)$ is regular if $D$ is regular and for any
$x \in M$ there is a topological ball $U(x) \subset M$ with $x \in U(M)$ and
 an orthonormal (with respect to the metric $g$)  family of
  smooth vector fields $\displaystyle \left\{ X_{1}, ..., X_{m}\right\}$ in
  $U(x)$ which satisfy the condition (C).
\label{dreg}
\end{definition}

The lenght of  a horizontal curve is
$$l(c) \ = \ \int_{a}^{b} \left(g_{c(t)} (\dot{c}(t),
\dot{c}(t))\right)^{\frac{1}{2}} \mbox{ d}t$$
The length depends on the metric $g$.

\begin{definition}
The Carnot-Carath\'eodory distance (or CC distance) associated to the sub-riemannian manifold is the
distance induced by the length $l$ of horizontal curves:
$$d(x,y) \ = \ \inf \left\{ l(c) \mbox{ : } c \in Hor(M,D) \
, \ c(a) = x \ , \  c(b) = y \right\} $$
\end{definition}

The Chow theorem ensures the existence of a horizontal path linking any two sufficiently
close points, therefore the CC distance is  locally finite. The distance
depends only on the distribution $D$ and metric $g$, and not on the choice of
vector fields $\displaystyle X_{1}, ..., X_{m}$ satisfying the condition (C).
The space $(M,d)$ is locally compact and complete, and the topology
induced by the distance $d$ is the same as the topology of the manifold $M$.
(These important details may be recovered from reading carefully the
constructive proofs of Chow theorem given by  Bella\"{\i}che
\cite{bell} or Gromov \cite{gromovsr}.)

\subsection{Normal frames}

Chow condition (C) is used to construct an adapted frame starting from a 
family of vector fields which generate the distribution $D$. A fundamental 
result in sub-riemannian geometry is the existence of normal frames. This 
existence result is based on the accumulation of various results by 
Bella\"{\i}che \cite{bell}, first to speak about normal frames, providing 
rigorous proofs for this existence in a flow of results between theorem 4.15 
and ending in the first half of section 7.3 (page 62), Gromov \cite{gromovsr} 
in his approximation theorem  p. 135 (conclusion of the point (a) below), as 
well in his convergence results concerning the nilpotentization of vector 
fields (related to point (b) below), Vodopyanov and others \cite{vodopis} 
\cite{vodopis2} \cite{vodokar} concerning the proof of basic results in 
sub-riemannian geometry under very weak regularity assumptions (for a 
discussion of this see \cite{buligasr}). There is no place here to submerge 
into this, we shall just assume that the object defined below exists.

In the following we stay in a small open neighbourhood of an arbitrary, but
fixed point $\displaystyle x_{0} \in M$. All results are local in
nature (that is they hold for some small open neighbourhood of an arbitrary, but
fixed point of the manifold $M$). That is why we shall no longer mention the
dependence of various objects on $\displaystyle x_{0}$, on the neighbourhood
$\displaystyle U(x_{0})$, or the distribution $D$.

We shall work further only with regular sub-riemannian manifolds, if not
otherwise stated. The topological dimension of $M$ is denoted by $n$, the
 step of the regular sub-riemannian manifold $(M,D,g)$ is denoted by $k$, the
dimension of the distribution is $m$, and there are numbers $\displaystyle
\nu_{j}$, $j = 1, ..., k$ such that for any $x \in M$ we have
$\displaystyle dim \, V^{j}(x) = \nu_{j}$. The Carnot-Carath\'eodory distance
is denoted by $d$.

An adapted frame $\displaystyle \left\{ X_{1}, ... , X_{n} \right\}$ is a
collection of smooth vector fields which is obtained by the construction
described below.

We start with a collection $X_{1}, ... , X_{m}$ of vector fields which satisfy
the condition (C). In particular  for any point $x$ the vectors  $\displaystyle
X_{1}(x), ... , X_{m}(x)$ form a basis for $\displaystyle D_{x}$.  We further
associate to any word $\displaystyle a_{1} .... a_{q}$ with letters in the
alphabet $1, ... ,m$ the  multi-bracket
$\displaystyle [X_{a_{1}}, [ ... , X_{a_{q}}] ... ]$.

One can add,  in the lexicographic order, $n-m$ elements to the set
$\displaystyle \left\{ X_{1}, ... , X_{m} \right\}$ until we get a collection
$\displaystyle \left\{ X_{1}, ... , X_{n} \right\}$ such that:
for any $j = 1, ..., k$ and for any point $x$ the set
$\displaystyle \left\{X_{1}(x), ..., X_{\nu_{j}}(x) \right\}$ is a basis for
$\displaystyle V^{j}(x)$.

Let $\displaystyle \left\{ X_{1}, ... , X_{n} \right\}$ be an adapted frame.
For any $j = 1, ..., n$  the degree $\displaystyle deg \, X_{j}$ of the vector
field $\displaystyle X_{j}$  is defined
as the only positive integer $p$ such that for any point $x$ we have
$$X_{j}(x) \in V^{p}_{x} \setminus V^{p-1}(x)$$

According 
with Gromov suggestions in the last section of Bella\"{\i}che \cite{bell}, 
the key details in the definition below are the uniform convergence assumptions.

\begin{definition}
An adapted frame $\displaystyle \left\{ X_{1}, ... , X_{n} \right\}$ is a {\bf normal
frame} if the following two conditions are satisfied:
\begin{enumerate}
\item[(a)] we have the limit
$$\lim_{\varepsilon \rightarrow 0_{+}} \frac{1}{\varepsilon} \, d\left(
\exp \left( \sum_{1}^{n}\varepsilon^{deg\, X_{i}} a_{i} X_{i} \right)(y), y \right) \ = \ A(y, a) \in
(0,+\infty)$$
which is uniform with respect to $y$ in compact sets and  vector $\displaystyle a=(a_{1}, ...,
a_{n}) \in W$, with $\displaystyle W \subset \mathbb{R}^{n}$ compact
neighbourhood of $\displaystyle 0 \in \mathbb{R}^{n}$,
\item[(b)] for any compact set $K\subset M$ with diameter (with respect to the
distance $d$) sufficiently small,  and for any $i = 1, ..., n$ there
are functions
$$ P_{i}(\cdot, \cdot, \cdot): U_{K} \times U_{K} \times K \rightarrow  \mathbb{R}$$
 with $\displaystyle U_{K} \subset \mathbb{R}^{n}$ a sufficiently small compact neighbourhood of
$\displaystyle 0 \in \mathbb{R}^{n}$ such that for any   $x \in K$ and any $\displaystyle
a,b \in U_{K}$ we have
$$\exp \left( \sum_{1}^{n} a_{i} X_{i} \right) (x) \ = \
\exp \left( \sum_{1}^{n} P_{i}(a, b, y) X_{i} \right) \circ \exp \left( \sum_{1}^{n}  b_{i} X_{i} \right) (x) $$
and such that the following limit exists
$$\lim_{\varepsilon \rightarrow 0_{+}}
\varepsilon^{-deg \, X_{i}} P_{i}(\varepsilon^{deg \, X_{j}} a_{j}, \varepsilon^{deg \, X_{k}} b_{k}, x)   \in
\mathbb{R}$$
and it is uniform with respect to $x  \in K$ and $\displaystyle a, b \in U_{K}$.
\end{enumerate}
\label{defnormal}
\end{definition}

In order to understand normal frames let us look to the
 case of a Lie group $G$ endowed with a left invariant distribution.
 The distribution is completely non-integrable if it is generated by the left
 translation of a vector subspace $D$ of the algebra
$\mathfrak{g} = T_{e}G$ which bracket generates the whole algebra
$\mathfrak{g}$. Take $\displaystyle \left\{ X_{1}, ..., X_{m}\right\}$ a collection of $m = \,
dim \, D$ left invariant independent vector fields and define with their help
an adapted frame. Then the adapted
frame $\displaystyle \left\{ X_{1}, ..., X_{n}\right\}$ is in fact normal.

With the help of a normal frame we can prove the existence of strong dilatation structures on regular sub-riemannian manifolds. The following is
a straightforward consequence of theorems 6.3, 6.4 \cite{buligasr}.

\begin{theorem}
 Let $(M,D, g)$ be a regular sub-riemannian manifold, $U \subset M$ an open set which admits a normal frame. Define for any $x \in  U$
and $\varepsilon > 0$ (sufficiently small if necessary),
the dilatation $\displaystyle \delta^{x}_{\varepsilon}$  given by:
$$\delta^{x}_{\varepsilon} \left(\exp\left( \sum_{i=1}^{n} a_{i} X_{i} \right)(x)\right) \  = \
\exp\left( \sum_{i=1}^{n} a_{i} \varepsilon^{deg X_{i}}  X_{i} \right)(x)$$
Then $(U, d,  \delta)$ is a strong dilatation structure.
\label{structhm}
\end{theorem}

\paragraph{Proof.}
Indeed, it is enough  to check the axioms A3, A4 of a strong dilatation structure, because the other axioms are obviously true.  By  theorem 6.3 \cite{buligasr} A3 is true and by theorem 6.4 \cite{buligasr} 
A4 is true. 
\hfill $\square$

\subsection{Carnot groups}

Carnot groups appear in sub-riemannian geometry
 as models of tangent spaces,   \cite{bell}, \cite{gromovsr}, \cite{pansu}. In particular such groups can be endowed with a structure
of sub-riemannian manifold.

Carnot groups are particular cases of normed conical groups.

\begin{definition}
A {\bf Carnot (or stratified homogeneous) group} $\displaystyle
(N, V_{1})$ is a pair  consisting of a real
connected simply connected group $N$  with  a distinguished subspace
$V_{1}$ of  the Lie algebra $Lie(N)$, such that  the following
direct sum decomposition occurs:
$$n \ = \ \sum_{i=1}^{m} V_{i} \ , \ \ V_{i+1} \ = \ [V_{1},V_{i}] $$
The number $m$ is the {\bf step of the group}. The number $\displaystyle Q \ = \ \sum_{i=1}^{m} i
\ dim V_{i}$ is called the {\bf homogeneous dimension of the group}.
\label{dccgroup}
\end{definition}

Because the group is nilpotent and simply connected, the
exponential mapping is a diffeomorphism. We shall identify the
group with the algebra, if is not locally otherwise stated.

The structure that we obtain is a set $N$ endowed with a Lie
bracket and a group multiplication operation, related by the
Baker-Campbell-Hausdorff formula. Remark that  the group operation
is polynomial.

Any Carnot group admits a one-parameter family of dilatations. For any
$\varepsilon > 0$, the associated dilatation is:
$$ x \ = \ \sum_{i=1}^{m} x_{i} \ \mapsto \ \delta_{\varepsilon} x \
= \ \sum_{i=1}^{m} \varepsilon^{i} x_{i} $$
Any such dilatation is a group morphism and a Lie algebra morphism.

In a Carnot group $N$ let us choose  an  euclidean norm
$\| \cdot \|$ on $\displaystyle V_{1}$.  We shall endow the group $N$ with
a structure of a sub-riemannian manifold. For this take the distribution
obtained from left translates of the space $V_{1}$. The metric on that
distribution is obtained by left translation of the inner product restricted to
$V_{1}$.

Because $V_{1}$  generates (the algebra) $N$
then any element $x \in N$ can be written as a product of
elements from $V_{1}$, in a controlled way, described in the following  useful lemma  (slight reformulation of Lemma 1.40,
Folland, Stein \cite{folstein}).

\begin{lemma}
Let $N$ be a Carnot group and $X_{1}, ..., X_{p}$ an orthonormal basis
for $V_{1}$. Then there is a
 a natural number $M$ and a function $g: \left\{ 1,...,M \right\}
\rightarrow \left\{ 1,...,p\right\}$ such that any element
$x \in N$ can be written as:
\begin{equation}
x \ = \ \prod_{i = 1}^{M} \exp(t_{i}X_{g(i)}) 
\label{fp2.4}
\end{equation}
Moreover, if $x$ is sufficiently close (in Euclidean norm) to
$0$ then each $t_{i}$ can be chosen such that $\mid t_{i}\mid \leq C
\| x \|^{1/m}$. 
\label{p2.4}
\end{lemma}

As a consequence we get:

\begin{corollary}
 The Carnot-Carath\'eodory distance
$$d(x,y) \ = \ \inf \left\{ \int_{0}^{1} \| c^{-1} \dot{c} \| \mbox{
d}t \ \mbox{ : } \ c(0) = x , \ c(1) = y , \quad \quad \right.$$
$$\left. \quad \quad \quad \quad  c^{-1}(t) \dot{c}(t) \in V_{1}
\mbox{ for a.e. } t \in [0,1]
\right\}$$
is finite for any two $x,y \in N$.  The distance is obviously left
invariant, thus it induces a norm on $N$.
\end{corollary}

The Carnot-Carath\'eodory distance induces a homogeneous norm on the Carnot group $N$ by the formula:
$\| x \| = d(0,x)$. From the invariance of the distance with respect to left translations we get: for any $x, y \in N$
$$\|x^{-1} y\| = d(x,y)$$

For any $x \in N$ and $\varepsilon > 0$ we define the dilatation $\displaystyle \delta^{x}_{\varepsilon} y = x \delta_{\varepsilon}(x^{-1} y)$.
Then $(N, d, \delta)$ is a dilatation structure, according to theorem \ref{tgrd}.

 Such dilatation structures have
the Radon-Nikodym property (defined further), as proven several times, in \cite{marmos1}, 
\cite{pansu}, or \cite{vodopis}.

\section{The Radon-Nikodym property}
\label{radon}

 Let $(X,d,\delta)$ be a strong dilatation structure or a length dilatation structure. We have then 
 a notion of differentiability for curves in $X$, seen as continuous functions from (a open interval 
 in) $\mathbb{R}$, with the usual dilatation structure, to $X$ with the dilatation structure $(X,d,\delta)$.  
 Further we want to see what  differentiability in the sense of definition \ref{defdiffer} means for curves. 
In proposition \ref{pintrinsicd} we shall arrive  to a kind of intrinsic notion of a distribution in a dilatation structure, with the geometrical meaning of a cone of all possible derivatives of curves passing through a point. 
 
\begin{definition}
In a normed conical group $N$   we shall denote by $D(N)$ the set of all 
$u\in N$ with the property that $\displaystyle 
\varepsilon \in ((0,\infty),+) \mapsto \delta_{\varepsilon} u \in N$ is a morphism of
groups.
\label{defdisn}
\end{definition}
$D(N)$ is always non empty, because it contains the neutral element of $N$. 
$D(N)$ is also a cone, with dilatations $\displaystyle   \delta_{\varepsilon}$, 
and a closed set.

 \begin{proposition}
 Let $(X,d,\delta)$ be a strong dilatation structure or a length dilatation structure and  
 let $\displaystyle c : [a, b] \rightarrow (X,d)$ be a continuous curve. For any 
 $x \in X$ and any $\displaystyle y \in T_{x}(X,d,\delta)$ we denote by 
 $$inv^{x}(y) \, = \, \Delta^{x}(y, x)$$ 
 the inverse of $y$ with respect to the group operation in  $\displaystyle  T_{x}(X,d,\delta)$. 
 Then the following are equivalent: 
 \begin{enumerate}
 \item[(a)] $c$ is derivable in $t \in (a,b)$ with respect to the dilatation structure $(X,d,\delta)$; 
 \item[(b)] there exists 
  $\displaystyle \dot{c}(t) \in D(T_{c(t}(X,d,\delta))$ such that 
$$\frac{1}{\varepsilon} d(c(t+\varepsilon) , \delta_{\varepsilon}^{c(t)} \dot{c}(t)) \rightarrow 0 $$
$$\frac{1}{\varepsilon} d(c(t-\varepsilon) , \delta_{\varepsilon}^{c(t)} 
inv^{c(t)}( \dot{c}(t))) \rightarrow 0  $$
\end{enumerate}
\label{pintrinsicd}
\end{proposition}

\paragraph{Proof.}  
 It is  straightforward that a conical group morphism $f: \mathbb{R} \rightarrow N$ is defined by its value $f(1)\in N$. Indeed, for any $a>0$ we have $\displaystyle f(a) = \delta_{a} f(1)$ and for any 
 $a<0$ we have $\displaystyle f(a) = \delta_{a} f(1)^{-1}$. From the morphism property we also 
 deduce that  
 $$\delta v = \left\{ \delta_{a} v \mbox{ : } a>0 , v=f(1) \mbox{ or } v=f(1)^{-1} \right\}$$
 is a one parameter group and that for all $\alpha, \beta >0$ we have 
 $\displaystyle \delta_{\alpha+\beta} u = \delta_{\alpha}u \,  \delta_{\beta}u$. We have
 therefore a bijection between conical group morphisms $f: \mathbb{R} \rightarrow
 (N,\delta)$  and elements of $D(N)$. 
 
 The curve  
 $\displaystyle c : [a,b] \rightarrow (X,d)$ is derivable in $t \in (a,b)$ if 
 and only if there is a morphism of normed conical groups 
 $\displaystyle f: \mathbb{R} \rightarrow T_{c(t}(X,d,\delta)$ such that 
 for any $a \in \mathbb{R}$ we have 
 $$\lim_{\varepsilon \rightarrow 0} \frac{1}{\varepsilon} \, d(c(t+\varepsilon a), 
 \delta^{c(t)}_{\varepsilon} f(a)) \ = \ 0$$
 Take  $\displaystyle \dot{c}(t) = f(1)$. Then $\displaystyle \dot{c}(t) \in 
 D(T_{c(t}(X,d,\delta))$. For any $a > 0$ we have $\displaystyle f(a) =
 \delta^{c(t)}_{a} \dot{c}(t)$; otherwise if $a < 0$ we have 
 $\displaystyle f(a) =
 \delta^{c(t)}_{a} \, inv^{c(t)} \,\dot{c}(t)$. This implies the equivalence 
 stated on the proposition.
 \hfill $\square$

 \begin{definition} 
A strong dilatation structure or a length dilatation structure
 $(X,d, \delta)$ has the {\bf Radon-Nikodym property (or rectifiability
 property, or RNP)} if  any  
Lipschitz curve $\displaystyle c : [a,b] \rightarrow (X,d)$ is derivable 
almost everywhere. 
\label{defrn}
 \end{definition}

\subsection{Two examples}

The following two easy examples will 
show that not any strong dilatation structure has the Radon-Nikodym property.

For  $\displaystyle (X,d)  =  ( \mathbb{V}, d)$, a real, finite dimensional,
normed vector space, with distance $d$ induced by the norm, the (usual) 
 dilatations $\displaystyle \delta^{x}_{\varepsilon}$ are given by:  
$$ \delta_{\varepsilon}^{x} y \ = \ x + \varepsilon (y-x) $$
Dilatations are defined everywhere.

There are few things to check:  axioms 0,1,2 are obviously 
true. For axiom A3, remark that for any $\varepsilon > 0$, $x,u,v \in X$ we 
have: 
$$\frac{1}{\varepsilon} d(\delta^{x}_{\varepsilon} u , 
\delta^{x}_{\varepsilon} v ) \ = \ d(u,v) \ , $$
therefore for any $x \in X$ we have $\displaystyle d^{x} = d$. 

Finally, let us check the axiom A4. For any $\varepsilon > 0$ and $x,u,v \in X$ we have
$$\delta_{\varepsilon^{-1}}^{\delta_{\varepsilon}^{x} u} \delta_{\varepsilon}^{x} v \ = \ 
x + \varepsilon  (u-x) + \frac{1}{\varepsilon} \left( x+ \varepsilon(v-x) - x - \varepsilon(u-x) \right) \ = \ $$
$$ = \ x + \varepsilon  (u-x) + v - u$$ 
therefore this quantity converges to 
$$x + v - u \ = \ x + (v - x) - (u - x)$$
as $\varepsilon \rightarrow 0$. The axiom A4 is verified. 

This dilatation structure has the Radon-Nikodym property.

Further is an example of a dilatation structure which does not have the
Radon-Nikodym property.  Take $\displaystyle X = \mathbb{R}^{2}$ with the euclidean 
distance $\displaystyle d$. For any $z \in \mathbb{C}$ of the 
form $z= 1+ i \theta$ we define dilatations 
$$\delta_{\varepsilon} x = \varepsilon^{z} x  \ .$$
It is easy to check that $\displaystyle (\mathbb{R}^{2},d, \delta)$ 
is a dilatation structure, with  dilatations 
$$\delta^{x}_{\varepsilon} y = x + \delta_{\varepsilon} (y-x)  $$

Two such dilatation structures (constructed with the help of complex numbers 
$1+ i \theta$ and $1+ i \theta'$) are equivalent if and only if $\theta = \theta'$.  

There are two other interesting  properties of these dilatation structures. 
The first is that if $\theta \not = 0$ then there are no non trivial 
Lipschitz curves in $X$ which are differentiable almost everywhere. It means
that such dilatation structure does not have the Radon-Nikodym property. 

The second property is that any holomorphic and Lipschitz function from $X$ to $X$ (holomorphic in the 
usual sense on $X = \mathbb{R}^{2} = \mathbb{C}$) is differentiable almost everywhere, but there are 
Lipschitz functions from $X$ to $X$ which are not differentiable almost everywhere (suffices to take a 
$\displaystyle \mathcal{C}^{\infty}$ function from  $\displaystyle \mathbb{R}^{2}$ to $\displaystyle \mathbb{R}^{2}$ which is not holomorphic).

\subsection{Length formula from Radon-Nikodym property}

 \begin{theorem}
Let $(X,d,\delta)$ be a strong dilatation structure with the Radon-Nikodym property, 
over a complete length metric space $(X,d)$. Then   for any $x, y \in X$ we have 
$$d(x,y) \ = \ \inf \left\{ \int_{a}^{b} d^{c(t)}(c(t),\dot{c}(t)) \mbox{ d}t  \mbox{ :
} c:[a,b]\rightarrow X \mbox{ Lipschitz }, \right. $$
$$\left.  c(a) = x , c(b) = y \right\}  $$
\label{fleng}
\end{theorem}

\paragraph{Proof.}
From theorem \ref{tupper} we deduce that for almost every $t\in(a,b)$ 
the upper dilatation of $c$ in $t$ can be expressed as:
$$Lip(c)(t) = \lim_{s\rightarrow t} \frac{d(c(s),c(t))}{\mid s-t \mid} $$
 
If the dilatation structure has the Radon-Nikodym property then for almost every $t \in [a,b]$ there is 
$\displaystyle \dot{c}(t) \in D(T_{c(t)} X)$ such that 
$$\displaystyle \frac{1}{\varepsilon} d(c(t+\varepsilon) , \delta_{\varepsilon}^{c(t)} \dot{c}(t)) 
\rightarrow 0$$ 
Therefore for almost every $t \in [a,b]$ we have
$$\displaystyle Lip(c)(t) = \lim_{\varepsilon\rightarrow 0} 
\frac{1}{\varepsilon} d(c(t+\varepsilon),c(t)) = d^{c(t)}(c(t),\dot{c}(t))$$ 
The formula for length follows from here. 
\hfill $\square$

A straightforward consequence is that the distance $d$ is uniquely determined by the 
"distribution" $\displaystyle x \in X \mapsto D(T_{x}(X,d,\delta))$ and the 
function which associates to any $x \in X$ the "norm" $\displaystyle 
\| \cdot \|_{x} : D(T_{x}(X,d,\delta)) \rightarrow [0,+ \infty)$. 
  
 \begin{corollary}
Let $(X,d,\delta)$ and $(X,\bar{d},\bar{\delta})$ be two strong dilatation structures 
with the Radon-Nikodym property , which are also complete length metric spaces, 
such that for any $x \in X$ we have  $\displaystyle D(T_{x}(X,d,\delta)) =
D(T_{x}(X,\bar{d}, \bar{\delta}))$ and $\displaystyle d^{x}(x,u ) = \bar{d}^{x}(x, u)$ for any 
$\displaystyle u \in  D(T_{x}(X,d,\delta))$. Then $d = \bar{d}$.
 \end{corollary}

\subsection{Equivalent dilatation structures and their distributions}

\begin{definition}
{\bf Two strong dilatation structures} $(X, \delta , d)$ and $(X,
\overline{\delta} , \overline{d})$  {\bf are equivalent}  if 
\begin{enumerate}
\item[(a)] the identity  map $\displaystyle id: (X, d) \rightarrow (X, \overline{d})$ is bilipschitz and 
\item[(b)]  for any $x \in X$ there are functions $\displaystyle P^{x}, Q^{x}$ (defined for $u \in X$ sufficiently close to $x$) such that  
\begin{equation}
\lim_{\varepsilon \rightarrow 0} \frac{1}{\varepsilon} \overline{d} \left( \delta^{x}_{\varepsilon} u ,  \overline{\delta}^{x}_{\varepsilon} Q^{x} (u) \right)  = 0 , 
\label{dequiva}
\end{equation}
\begin{equation}
 \lim_{\varepsilon \rightarrow 0} \frac{1}{\varepsilon} d \left( \overline{\delta}^{x}_{\varepsilon} u ,  
 \delta^{x}_{\varepsilon} P^{x} (u) \right)  = 0 , 
\label{dequivb}
\end{equation}
uniformly with respect to $x$, $u$ in compact sets. 
\end{enumerate}
\label{dilequi}
\end{definition}

\begin{proposition}
 $(X, \delta , d)$ and $(X, \overline{\delta} , \overline{d})$  are equivalent  if and 
only if 
\begin{enumerate}
\item[(a)] the identity  map $\displaystyle id: (X, d) \rightarrow (X,
\overline{d})$ is bilipschitz, 
\item[(b)]  for any $x \in X$ there are conical group morphisms: 
 $$\displaystyle P^{x}: T_{x}(X, \overline{\delta} , \overline{d}) 
 \rightarrow T_{x} (X, \delta , d) \mbox{ and } \displaystyle  Q^{x}: T_{x} (X, \delta , d) \rightarrow 
 T_{x}(X, \overline{\delta} , \overline{d})$$
  such that the following limits exist  
\begin{equation}
\lim_{\varepsilon \rightarrow 0}  \left(\overline{\delta}^{x}_{\varepsilon}\right)^{-1}  \delta^{x}_{\varepsilon} (u) = Q^{x}(u) , 
\label{dequivap}
\end{equation}
\begin{equation}
 \lim_{\varepsilon \rightarrow 0}  \left(\delta^{x}_{\varepsilon}\right)^{-1}  \overline{\delta}^{x}_{\varepsilon} (u) = P^{x}(u) , 
\label{dequivbp}
\end{equation}
and are uniform with respect to $x$, $u$ in compact sets. 
\end{enumerate}
\label{pdilequi}
\end{proposition}

 The next theorem shows a link between the tangent bundles of equivalent dilatation structures. 
 
 \begin{theorem} 
 Let $(X, d, \delta)$ and $(X, \overline{d}, \overline{\delta})$  be  equivalent
 strong  dilatation structures.  Then for any $x \in X$ and 
 any $u,v \in X$ sufficiently close to $x$ we have:
 \begin{equation}
 \overline{\Sigma}^{x}(u,v) = Q^{x} \left( \Sigma^{x} \left( P^{x}(u) , P^{x}(v) \right)\right) . 
 \label{isoequiv}
 \end{equation}
 The two tangent bundles  are therefore isomorphic in a natural sense. 
 \label{tisoequiv}
 \end{theorem}
 As a consequence, the following corollary is straightforward. 
 
\begin{corollary}
Let  $(X, d, \delta)$ and $(X, \overline{d}, \overline{\delta})$  be  
equivalent strong dilatation structures. Then for any $x \in X$ we have 
$$\displaystyle Q^{x} (D(T_{x}(X, \delta , d))) \ = \ D(T_{x}(X, \overline{\delta} ,
\overline{d}))  $$

If  $(X, d, \delta)$ has the Radon-Nikodym property , then 
$(X, \overline{d}, \overline{\delta})$ has the same property. 

Suppose that $(X, d, \delta)$ and $(X, \overline{d}, \overline{\delta})$  are 
complete length spaces with the Radon-Nikodym property . If the functions 
$\displaystyle P^{x}, Q^{x}$ from definition \ref{dilequi} (b) are isometries, 
then $\displaystyle d = \overline{d}$. 
\end{corollary}

 \section{Tempered dilatation structures}
\label{stemp}

The notion of a tempered dilatation structure is inspired by 
the results from Venturini \cite{venturini} and Buttazzo, De Pascale and 
Fragala \cite{buttazzo1}. 

The examples of length dilatation structures from this section are provided 
by the extension of some results from \cite{buttazzo1} (propositions 2.3, 2.6 
and a part of theorem 3.1) to dilatation 
structures. 

The following definition gives a class of distances $\mathcal{D}(\Omega,
\bar{d}, \bar{\delta})$, associated to a strong dilatation structure 
$(\Omega, \bar{d}, \bar{\delta})$, which generalizes the class of distances 
$\mathcal{D}(\Omega)$ from \cite{buttazzo1}, definition 2.1. 

\begin{definition}
For any strong dilatation structure $(\Omega, \bar{d}, \bar{\delta})$   
 we define the class $\mathcal{D}(\Omega,
\bar{d}, \bar{\delta})$ of all distance functions $d$ on $\Omega$ such that 
\begin{enumerate}
\item[(a)] $d$ is a length distance, 
\item[(b)] for any $\varepsilon > 0$ and any $x, u, v$ sufficiently 
close the are constants $0 < c < C$ such that: 
\begin{equation}
c \, \bar{d}^{x}(u,v) \, \leq \, \frac{1}{\varepsilon} \,
d(\bar{\delta}^{x}_{\varepsilon} u , \bar{\delta}^{x}_{\varepsilon} v ) \, \leq 
\, C \, \bar{d}^{x}(u,v) 
\label{new2.3}
\end{equation}
\end{enumerate}
The dilatation structure $(\Omega, \bar{d}, \bar{\delta})$ is {\bf tempered} if 
 $\bar{d} \in \mathcal{D}(\Omega,
\bar{d}, \bar{\delta})$. 

On $\mathcal{D}(\Omega, \bar{d}, \bar{\delta})$ we put the topology of uniform
convergence (induced by distance $\bar{d}$) on compact subsets of $\Omega \times
\Omega$. 
\label{dtempered}
\end{definition}

To any distance $d \in \mathcal{D}(\Omega, \bar{d}, \bar{\delta})$ we associate
the function: 
$$\phi_{d}(x, u) \, = \, \limsup_{\varepsilon \rightarrow 0}
\frac{1}{\varepsilon} \, d(x, \delta^{x}_{\varepsilon} u ) $$
defined for any $x, u \in \Omega$ sufficiently close. We have therefore 
\begin{equation}
c \, \bar{d}^{x}(x, u) \, \leq \, \phi_{d}(x,u) \, \leq \, C \, \bar{d}^{x}(x,u)
\label{new2.6}
\end{equation}

Notice that if $d \in \mathcal{D}(\Omega, \bar{d}, \bar{\delta})$ then for any 
$x, u, v$ sufficiently close we have 
$$- \bar{d}(x,u) \, O(\bar{d}(x,u)) \, + \, c \,  \bar{d}^{x}(u,v) \, \leq $$
$$ \leq 
\, d(u,v) \, \leq \,   \, C \, \bar{d}^{x}(u,v) \, + \, 
 \bar{d}(x,u) \, O(\bar{d}(x,u))$$

If $c: [0,1] \rightarrow \Omega$ is a $d$-Lipschitz curve and 
$d \in \mathcal{D}(\Omega, \bar{d}, \bar{\delta})$ then we may decompose it 
in a finite family of curves $\displaystyle c_{1}, ... , c_{n}$ (with $n$ depending on $c$) 
such that there are $\displaystyle x_{1}, ... , x_{n} \in \Omega$ with 
$\displaystyle c_{k} $ is $\displaystyle \bar{d}^{x_{k}}$-Lipschitz. Indeed, the image of the 
curve $c([0,1])$ is compact, therefore we may cover it with a finite number of
balls $\displaystyle B(c(t_{k}), \rho_{k}, \bar{d}^{c(t_{k})})$ and apply 
(\ref{new2.3}). If moreover $(\Omega, \bar{d}, \bar{\delta})$ is tempered then 
it follows that $c: [0,1] \rightarrow \Omega$  $d$-Lipschitz curve is equivalent
with $c$ $\bar{d}$-Lipschitz curve. 

By using the same arguments as in the proof of theorem \ref{fleng}, we get the
following extension of proposition 2.4 \cite{buttazzo1}. 

\begin{proposition}
If $(\Omega, \bar{d}, \bar{\delta})$ is tempered, with the Radon-Nikodym 
property, and $d \in \mathcal{D}(\Omega, \bar{d}, \bar{\delta})$ then 
$$d(x,y) \ = \ \inf \left\{ \int_{a}^{b} \phi_{d}(c(t),\dot{c}(t)) \mbox{ d}t  \mbox{ :
} c:[a,b]\rightarrow X \mbox{ $\bar{d}$-Lipschitz }, \right. $$
$$\left.  c(a) = x , c(b) = y \right\}  $$
\label{new2.4}
\end{proposition}

The next theorem is a generalization of a part of theorem 3.1
\cite{buttazzo1}.  

\begin{theorem}
Let $(\Omega, \bar{d}, \bar{\delta})$ be a strong dilatation structure which is 
tempered, with the Radon-Nikodym 
property, and $\displaystyle d_{n} \in \mathcal{D}(\Omega, \bar{d}, \bar{\delta})$ 
 a sequence  of distances converging to $d \in \mathcal{D}(\Omega, \bar{d},
 \bar{\delta})$. Denote by $\displaystyle L_{n}, L$ the length functional induced
 by the distance $\displaystyle d_{n}$, respectively by $d$. 
 Then $\displaystyle L_{n}$ $\Gamma$-converges to $L$. 
 \label{new3.1}
 \end{theorem}
 
 \paragraph{Proof.}
 This is the generalization of the implication (i) $\Rightarrow$ (iii), theorem 
 3.1 \cite{buttazzo1}. The proof (p. 252-253) is almost identical, we only need to 
replace everywhere expressions like  $\mid x - y\mid$ by $\bar{d}(x,y)$ and 
use proposition \ref{new2.4}, relations (\ref{new2.6}) and (\ref{new2.3})
instead of respectively proposition 2.4 and relations (2.6) and (2.3)
\cite{buttazzo1}. 
 \hfill $\square$

Using this result we obtain a large class of examples of length dilatation
structures.

\begin{corollary}
If  $(\Omega, \bar{d}, \bar{\delta})$ is strong dilatation structure which is 
tempered and it has  the Radon-Nikodym 
property then it is a length dilatation structure. 
\label{cortemp}
\end{corollary}

\paragraph{Proof.}
Indeed, from the hypothesis we deduce that $\displaystyle
\bar{\delta}^{x}_{\varepsilon} \bar{d} \, \in \, 
 \mathcal{D}(\Omega, \bar{d}, \bar{\delta})$. For any sequence 
 $\displaystyle \varepsilon_{n} \rightarrow 0$ we thus obtain a sequence of 
 distances $\displaystyle d_{n} \, = \, \bar{\delta}^{x}_{\varepsilon_{n}} 
 \bar{d}$ converging to $\bar{d}^{x}$. We apply now theorem \ref{new3.1} 
 and we get the result. 
\hfill $\square$

\section{Coherent projections}
\label{cohp}

For a given dilatation structure with the Radon-Nikodym property, we shall give a procedure to 
construct another dilatation structure, such that the first one looks down to 
the the second one. 

This will be done with the help of coherent projections.

\begin{definition}
Let $(X,\bar{d}, \bar{\delta})$ be a strong dilatation structure. A {\bf
coherent projection}  
of $(X,\bar{d}, \bar{\delta})$  is a function which associates to any $x \in X$ and $\varepsilon
\in (0,1]$ a map 
$\displaystyle Q^{x}_{\varepsilon} : U(x) \rightarrow X$ such that: 
\begin{enumerate}
\item[(I)] $\displaystyle Q^{x}_{\varepsilon} : U(x) \rightarrow  Q^{x}_{\varepsilon}(U(x))$ 
is invertible and the inverse will be denoted by $\displaystyle Q^{x}_{\varepsilon^{-1}}$; for
any $\varepsilon, \mu > 0$ and any $x \in X$ we have
$$Q_{\varepsilon}^{x} \, \bar{\delta}^{x}_{\mu} \ = \ \bar{\delta}^{x}_{\mu} \,
Q^{x}_{\varepsilon}$$
\item[(II)] the limit $\displaystyle \lim_{\varepsilon \rightarrow 0} 
Q^{x}_{\varepsilon} u \ = \ Q^{x} u$  
is uniform with respect to $x, u$ in compact sets. 
\item[(III)] for any $\varepsilon, \mu > 0$ and any $x \in X$ we have $\displaystyle 
Q^{x}_{\varepsilon} \, Q^{x}_{\mu} \ = \ Q^{x}_{\varepsilon \mu}$. Also $\displaystyle
Q^{x}_{1} \ = \ id$ and $\displaystyle Q^{x}_{\varepsilon} x \ = \ x$.

\item[(IV)] define $\displaystyle \Theta_{\varepsilon}^{x}(u,v) \ = \ 
\bar{\delta}^{x}_{\varepsilon^{-1}} \,
Q^{\bar{\delta}^{x}_{\varepsilon} Q^{x}_{\varepsilon} u}_{\varepsilon^{-1}} \,
\bar{\delta}^{x}_{\varepsilon} Q^{x}_{\varepsilon} v$.  
Then the limit exists
$$\lim_{\varepsilon \rightarrow 0} \Theta^{x}_{\varepsilon}(u,v) \ = \ \Theta^{x}(u,v)$$
and it is uniform with respect to $x, u , v$ in compact sets.
\end{enumerate}
\label{defcoh}
\end{definition}

\begin{remark}
Property (IV) is basically a smoothness condition on the coherent projection
$Q$, relative to the strong dilatation structure $(X, \bar{d}, \bar{\delta})$.
\end{remark}

\begin{proposition}

Let $(X,\bar{d}, \bar{\delta})$ be a strong dilatation structure and $Q$ a coherent 
projection. We define $\displaystyle \delta^{x}_{\varepsilon} \ = \
\bar{\delta}^{x}_{\varepsilon} \, Q^{x}_{\varepsilon}$.  Then: 
\begin{enumerate}
\item[(a)]  for any $\varepsilon, \mu > 0$ and any $x \in X$ we have $\displaystyle 
\delta^{x}_{\varepsilon} \, \bar{\delta}^{x}_{\mu} \ = \  \bar{\delta}^{x}_{\mu} \, 
\delta^{x}_{\varepsilon}$. 
\item[(b)] for any $x \in X$ we have $\displaystyle Q^{x} \, Q^{x} \ = \ Q^{x}$ (thus 
$\displaystyle Q^{x}$ is a projection). 
\item[(c)]  $\delta$ satisfies the 
conditions A1, A2, A4 from definition \ref{defweakstrong}. 
\end{enumerate}
\label{p1proj}
\end{proposition}

\paragraph{Proof.}
(a) this is a consequence of the commutativity condition (I) (second part). Indeed, we have 
$\displaystyle \delta^{x}_{\varepsilon} \, \bar{\delta}^{x}_{\mu} \ = \ 
\bar{\delta}^{x}_{\varepsilon} \, Q^{x}_{\varepsilon} \, \bar{\delta}^{x}_{\mu} \ = \ 
\bar{\delta}^{x}_{\varepsilon} \, \bar{\delta}^{x}_{\mu} \, Q^{x}_{\varepsilon} \ = \ 
\bar{\delta}^{x}_{\mu} \, \bar{\delta}^{x}_{\varepsilon} \, Q^{x}_{\varepsilon} \ = \ 
 \bar{\delta}^{x}_{\mu} \, \delta^{x}_{\varepsilon}$. 

(b) we pass to the limit $\varepsilon \rightarrow 0$  in 
 the equality $\displaystyle Q^{x}_{\varepsilon^{2}} \ = \ Q^{x}_{\varepsilon}
\, Q^{x}_{\varepsilon}$ and we get, based on condition (II), that $\displaystyle Q^{x} \, Q^{x}
\ = \ Q^{x}$.

(c) Axiom A1 for $\delta$ is equivalent with (III). Indeed,  the equality 
$\displaystyle \delta^{x}_{\varepsilon} \, \delta^{x}_{\mu} \, = \, 
\delta^{x}_{\varepsilon \mu}$  is 
equivalent with: 
$\displaystyle \bar{\delta}^{x}_{\varepsilon \mu} \, Q^{x}_{\varepsilon \mu} \, = \, 
  \bar{\delta}^{x}_{\varepsilon \mu} \, Q^{x}_{\varepsilon} \, Q^{x}_{\mu}$. 
  This is true because $\displaystyle 
Q^{x}_{\varepsilon} \, Q^{x}_{\mu} \, = \, Q^{x}_{\varepsilon \mu}$. 
We also have $\displaystyle  \delta^{x}_{1} \, =  \, \bar{\delta}^{x}_{1} Q^{x}_{1} 
\, = \, Q^{x}_{1} \, = \, id $. 
Moreover $\displaystyle \delta^{x}_{\varepsilon} x \, = \, 
\bar{\delta}^{x}_{\varepsilon} \, Q^{x}_{\varepsilon} x \, = \, 
Q^{x}_{\varepsilon} \,  \bar{\delta}^{x}_{\varepsilon} x \, = \, 
Q^{x}_{\varepsilon} x  \, = \, x$. 
Let us compute now: 
$$\Delta^{x}_{\varepsilon} (u,v) \ = \ \delta^{\delta^{x}_{\varepsilon} u}_{\varepsilon^{-1}} 
\, \delta^{x}_{\varepsilon} v \ = \ \bar{\delta}^{\delta^{x}_{\varepsilon} u}_{\varepsilon^{-1}} 
\, Q^{\delta^{x}_{\varepsilon} u}_{\varepsilon^{-1}} \, \delta^{x}_{\varepsilon} v \ =
\ $$ 
$$ = \ \bar{\delta}^{\delta^{x}_{\varepsilon} u}_{\varepsilon^{-1}} \, \bar{\delta}^{x}_{\varepsilon} 
\, \Theta_{\varepsilon}^{x}(u,v) \ = \ \bar{\Delta}^{x}_{\varepsilon}(Q^{x}_{\varepsilon} u, 
\Theta_{\varepsilon}^{x}(u,v))$$
We can pass to the limit in the last term of this string of equalities and 
we prove that   the axiom A4 is satisfied by $\delta$: there exists the limit 
\begin{equation}
\Delta^{x}(u,v) \, = \, \lim_{\varepsilon \rightarrow 0} \Delta_{\varepsilon}^{x}(u,v) 
\label{defstu}
\end{equation}
which is uniform as written in A4, 
  moreover  we have the equality 
\begin{equation}
\Theta_{\varepsilon}^{x}(u,v) \ = \ \bar{\Sigma}^{x}_{\varepsilon}(Q^{x}_{\varepsilon} u, 
\Delta^{x}_{\varepsilon} (u,v))  
\label{neednext1}
\end{equation}
\hfill $\square$

We collect two useful relations in the next proposition.

\begin{proposition}
Let $(X,\bar{d}, \bar{\delta})$ be a strong dilatation structure and $Q$ a coherent 
projection. We denote by $\delta$ the field of dilatations induced by the coherent projection, as in the previous proposition, and by  $\Delta^{x}$ is defined by (\ref{defstu}). Then we have:  
\begin{equation}
\Delta^{x} (u,v) \, = \, \bar{\Delta}^{x}(Q^{x} u, 
\Theta^{x}(u,v)) 
\label{recon}
\end{equation}
\begin{equation}
Q^{x} \Delta^{x}(u,v) \, = \, \bar{\Delta}^{x}(Q^{x} u, Q^{x} v)
\label{morph}
\end{equation}
\end{proposition}

\paragraph{Proof.}
After passing to the limit with $\varepsilon \rightarrow 0$ in the  relation
(\ref{neednext1})  we get the formula (\ref{recon}). In order to prove 
(\ref{morph}) we notice that: 
$$Q^{\delta^{x}_{\varepsilon} u}_{\varepsilon} \, 
\Delta^{x}_{\varepsilon}(u,v) \, = \, 
 Q^{\delta^{x}_{\varepsilon} u}_{\varepsilon} \delta^{\delta^{x}_{\varepsilon}
 u}_{\varepsilon^{-1}} \delta^{x}_{\varepsilon} v \, =  $$ 
 $$ = \,  \bar{\delta}^{\delta^{x}_{\varepsilon} u}_{\varepsilon^{-1}}
 \bar{\delta}^{x}_{\varepsilon}  
 Q^{x}_{\varepsilon} v \, = \,
 \bar{\Delta}^{x}_{\varepsilon} (Q^{x}_{\varepsilon} u,
 Q^{x}_{\varepsilon} v)$$ 
which gives(\ref{morph}) as we 
 pass to the limit with $\varepsilon \rightarrow 0$ in this relation. 
 \hfill $\square$

Next is described the notion of $Q$-horizontal curve.

\begin{definition}
Let $(X,\bar{d}, \bar{\delta})$ be a strong dilatation structure and $Q$ a 
coherent projection. A curve $c: [a,b] \rightarrow X$ is $Q$-{\bf horizontal} if 
for almost any $t \in [a,b]$ the curve $c$ is derivable and the derivative 
of $c$ at $t$, denoted by $\dot{c}(t)$ has the property:
\begin{equation}
 Q^{c(t)} \dot{c}(t) \, = \, \dot{c}(t)
 \label{horprop}
 \end{equation}
 A curve $c: [a,b] \rightarrow X$ is $Q$-{\bf everywhere horizontal} if for all 
 $t \in [a,b]$ the curve $c$ is derivable and the derivative has the
 horizontality property (\ref{horprop}).
\end{definition}

We shall now use the notations from section \ref{seclds}. We  look first at
some induced dilatation structures. 

For any $x \in X$ and $\varepsilon \in (0,1)$ the  dilatation 
$\displaystyle \delta^{x}_{\varepsilon}$ can be seen as an isomorphism of 
strong dilatation structures with coherent projections: 
$$\delta^{x}_{\varepsilon} : (U(x), \delta^{x}_{\varepsilon} \bar{d} ,
\hat{\delta}^{x}_{\varepsilon}, \hat{Q}^{x}_{\varepsilon}) \rightarrow 
(\delta^{x}_{\varepsilon} U(x) , \frac{1}{\varepsilon} \bar{d} , \bar{\delta}, Q)$$
which defines the dilatations $\displaystyle \hat{\delta}^{x,
\cdot}_{\varepsilon, \cdot}$ and
coherent projection $\displaystyle \hat{Q}^{x}_{\varepsilon}$ by: 

$$\hat{\delta}^{x, u}_{\varepsilon, \mu} \ = \ \delta^{x}_{\varepsilon^{-1}} \, 
\bar{\delta}_{\mu}^{\delta^{x}_{\varepsilon} u} \, \delta^{x}_{\varepsilon} 
$$ 
$$\hat{Q}^{x, u}_{\varepsilon, \mu} \ = \ \delta^{x}_{\varepsilon^{-1}} \, 
Q_{\mu}^{\delta^{x}_{\varepsilon} u} \, \delta^{x}_{\varepsilon} 
 $$
Also the dilatation $\displaystyle \bar{\delta}^{x}_{\varepsilon}$
is an isomorphism of strong dilatation structures with coherent projections: 
$$\bar{\delta}^{x}_{\varepsilon} : (U(x), \bar{\delta}^{x}_{\varepsilon} \bar{d},
\bar{\delta}^{x}_{\varepsilon}, \bar{Q}^{x}_{\varepsilon}) \rightarrow 
(\bar{\delta}^{x}_{\varepsilon} U(x) , \frac{1}{\varepsilon} \bar{d} , 
\bar{\delta}, Q)$$
which defines the dilatations $\displaystyle \bar{\delta}^{x,
\cdot}_{\varepsilon, \cdot}$ and
coherent projection $\displaystyle \bar{Q}^{x}_{\varepsilon}$ by: 
$$\bar{\delta}^{x, u}_{\varepsilon, \mu} \ = \ \bar{\delta}^{x}_{\varepsilon^{-1}} \, 
\bar{\delta}_{\mu}^{\bar{\delta}^{x}_{\varepsilon} u} \, \bar{\delta}^{x}_{\varepsilon} 
 $$
$$\bar{Q}^{x, u}_{\varepsilon, \mu} \ = \ \bar{\delta}^{x}_{\varepsilon^{-1}} \, 
Q_{\mu}^{\bar{\delta}^{x}_{\varepsilon} u} \, \bar{\delta}^{x}_{\varepsilon} 
 $$
Because $\displaystyle \delta^{x}_{\varepsilon} \, = \,
\bar{\delta}^{x}_{\varepsilon} \, Q^{x}_{\varepsilon}$ we get that 
$$Q^{x}_{\varepsilon}: (U(x), \delta^{x}_{\varepsilon} \bar{d} ,
\hat{\delta}^{x}_{\varepsilon}, \hat{Q}^{x}_{\varepsilon}) \rightarrow 
(Q^{x}_{\varepsilon} U(x), \bar{\delta}^{x}_{\varepsilon} d ,
\bar{\delta}^{x}_{\varepsilon}, \bar{Q}^{x}_{\varepsilon})$$ 
is an isomorphism of strong dilatation structures with coherent projections.

Further is a useful description of the coherent projection 
$\displaystyle \hat{Q}^{x}_{\varepsilon}$. 

\begin{proposition}
With the  notations previously made, for any $\varepsilon \in (0,1]$, $x, u, v \in
X$ sufficiently close and $\mu > 0$ we have:
\begin{enumerate}
\item[(i)] $\displaystyle \hat{Q}^{x, u}_{\varepsilon, \mu} v \, = \, 
\Sigma^{x}_{\varepsilon}(u, Q_{\mu}^{\delta^{x}_{\varepsilon} u}
\Delta^{x}_{\varepsilon}(u,v))$, 
\item[(ii)] $\displaystyle \hat{Q}^{x, u}_{\varepsilon} v \, = \, 
\Sigma^{x}_{\varepsilon}(u, Q^{\delta^{x}_{\varepsilon} u}
\Delta^{x}_{\varepsilon}(u,v))$. 
\end{enumerate}
\label{pdesc}
\end{proposition}

\paragraph{Proof.}
(i) implies (ii) when $\mu \rightarrow 0$, thus it is sufficient to prove only 
the first point. This is the result of a computation: 
$$\hat{Q}^{x, u}_{\varepsilon, \mu} v \, = \, \delta^{x}_{\varepsilon^{-1}} \, 
Q_{\mu}^{\delta^{x}_{\varepsilon} u} \, \delta^{x}_{\varepsilon} \, = \, $$
$$= \, \delta^{x}_{\varepsilon^{-1}} \,
\delta_{\varepsilon}^{\delta^{x}_{\varepsilon} u} \, 
Q_{\mu}^{\delta^{x}_{\varepsilon} u} \, 
\delta_{\varepsilon^{-1}}^{\delta^{x}_{\varepsilon} u} \,\delta^{x}_{\varepsilon} \, = \, 
\Sigma^{x}_{\varepsilon}(u, Q_{\mu}^{\delta^{x}_{\varepsilon} u}
\Delta^{x}_{\varepsilon}(u,v)) $$
\hfill $\square$

{\bf Notation concerning derivatives.} We shall denote the derivative of a curve with respect to the dilatations 
$\displaystyle \hat{\delta}^{x}_{\varepsilon}$ by $\displaystyle 
\frac{\hat{d}^{x}_{\varepsilon}}{dt}$. Also, the derivative 
of the curve $c$ with respect to $\bar{\delta}$ is denoted by 
$\displaystyle \frac{\bar{d}}{dt}$, and so on.

 By computation we get: the curve $c$ is 
$\displaystyle \hat{\delta}^{x}_{\varepsilon}$-derivable if and only if 
$\displaystyle \delta^{x}_{\varepsilon} c$ is $\bar{\delta}$-derivable and 
$$ \frac{\hat{d}^{x}_{\varepsilon}}{dt} \, c(t) \ = \
\delta^{x}_{\varepsilon^{-1}} \, \frac{\bar{d}}{dt} \left(
\delta^{x}_{\varepsilon} c \right) (t) $$

With these notations we give a proposition which explains that the operator 
$\displaystyle \Theta^{x}_{\varepsilon}$ ,
from the definition of coherent projections, is a lifting operator. 

\begin{proposition}
 If the curve $\displaystyle \delta^{x}_{\varepsilon} c$ is $Q$-horizontal then 
 $$\frac{\bar{d}^{x}_{\varepsilon}}{dt} \left( Q^{x}_{\varepsilon} c \right) (t) \, =
 \, \Theta^{x}_{\varepsilon}(c(t), \frac{\hat{d}^{x}_{\varepsilon}}{dt}  c  (t))
 $$
 \end{proposition}

 \paragraph{Proof.}If the curve 
 $\displaystyle Q^{x}_{\varepsilon} c$ is $\displaystyle
 \bar{\delta}^{x}_{\varepsilon}$ derivable and $\displaystyle
 \bar{Q}^{x}_{\varepsilon}$ horizontal. We have therefore: 
 $$\frac{\bar{d}^{x}_{\varepsilon}}{dt} \left( Q^{x}_{\varepsilon} c \right) (t) 
 \, = \, \bar{\delta}^{x}_{\varepsilon^{-1}} \, Q^{\delta^{x}_{\varepsilon} c(t)} \,
 \bar{\delta}^{x}_{\varepsilon} \, \frac{\bar{d}^{x}_{\varepsilon}}{dt} \left(
 Q^{x}_{\varepsilon} c \right) (t)$$
 which implies: 
 $$\bar{\delta}^{x}_{\varepsilon} \, 
 \frac{\bar{d}^{x}_{\varepsilon}}{dt} \left( Q^{x}_{\varepsilon} c \right) (t) \, =
 \, Q^{\delta^{x}_{\varepsilon} c(t)}_{\varepsilon^{-1}} \, \bar{\delta}^{x}_{\varepsilon} \, \frac{\bar{d}^{x}_{\varepsilon}}{dt} \left(
 Q^{x}_{\varepsilon} c \right) (t) \, = \, 
 Q^{\delta^{x}_{\varepsilon} c(t)}_{\varepsilon^{-1}} \, \delta^{x}_{\varepsilon} \,
 \frac{\hat{d}^{x}_{\varepsilon}}{dt}  c  (t)$$ 
 which is the formula we wanted to prove. 
 \hfill $\square$

\subsection{Distributions in sub-riemannian spaces}

The inspiration for the notion of coherent projection comes from sub-riemannian
geometry. We shall look to the section \ref{srsmooth} with a fresh eye. 

Further we shall work locally, just as in the mentioned section. Same notations
are used.  
Let $\displaystyle \left\{ Y_{1}, ..., Y_{n} \right\}$ be a frame induced by 
a parameterization $\displaystyle \phi: O \subset \mathbb{R}^{n} \rightarrow U
\subset M$ of a small open, connected set $U$ in the manifold $M$. This
parameterization induces a affine dilatation structure on $U$, by
$$\tilde{\delta}^{\phi(a)}_{\varepsilon} \, \phi(b) \ = \ \phi\left( a +
\varepsilon(-a+b)\right) $$
We take the distance $\tilde{d}(\phi(a), \phi(b)) \, = \, \|b-a\|$. 

Let $\displaystyle \left\{ X_{1}, ..., X_{n} \right\}$ be a normal frame, 
cf. definition \ref{defnormal}, $d$ be the Carnot-Carath\'eodory distance 
and 
$$\delta^{x}_{\varepsilon} \left(\exp\left( \sum_{i=1}^{n} a_{i} X_{i} \right)(x)\right) \  = \
\exp\left( \sum_{i=1}^{n} a_{i} \varepsilon^{deg X_{i}}  X_{i} \right)(x)$$
be the dilatation structure associated, cf. theorem \ref{structhm}. 

We may take another dilatation structure, constructed as follows: extend the
metric $g$ on the distribution $D$ to a riemannian metric on $M$, denoted 
for convenience also by $g$. Let $\displaystyle \bar{d}$ be the riemannian
distance induced by the riemannian metric $g$, and the dilatations 
$$\bar{\delta}^{x}_{\varepsilon} \left(\exp\left( \sum_{i=1}^{n} a_{i} X_{i} \right)(x)\right) \  = \
\exp\left( \sum_{i=1}^{n} a_{i} \varepsilon  X_{i} \right)(x)$$
Then $\displaystyle (U, \bar{d}, \bar{\delta})$ is a strong dilatation 
structure which is equivalent with the dilatation structure 
$\displaystyle (U, \tilde{d}, \tilde{\delta})$. 

From now we may define coherent projections associated either to the pair 
$\displaystyle (\tilde{\delta}, \delta)$ or to the pair $\displaystyle 
(\bar{\delta}, \delta)$. Because we put everything on the manifold (by 
the use of the chosen frames), we shall obtain different coherent projections, 
both inducing the same dilatation structure $(U,d, \delta)$. 

 Let us define $\displaystyle Q^{x}_{\varepsilon}$ by: 
\begin{equation}
Q^{x}_{\varepsilon}\,  \left(\exp\left( \sum_{i=1}^{n} a_{i} X_{i} \right)(x)\right) \  = \
\exp\left( \sum_{i=1}^{n} a_{i} \varepsilon^{deg X_{i} - 1}  X_{i} \right)(x)
\label{eqsr1}
\end{equation}

\begin{proposition}
$Q$ is a coherent projection associated with the dilatation structure
$\displaystyle (U, \bar{d}, \bar{\delta})$ .
\label{exsr}
\end{proposition}

\paragraph{Proof.}(I) definition \ref{defcoh} is true, because 
$\displaystyle \delta^{x}_{\varepsilon} \, u \, = \, Q^{x}_{\varepsilon} \, 
\bar{\delta}^{x}_{\varepsilon}$ and $\displaystyle \delta^{x}_{\varepsilon} \, 
\bar{\delta}^{x}_{\varepsilon} \, = \, \bar{\delta}^{x}_{\varepsilon}
\delta^{x}_{\varepsilon}$. (II), (III) and (IV) are consequences of these facts 
and theorem \ref{structhm}, with a proof similar to the one of proposition 
\ref{p1proj}. 
\hfill $\square$

Definition (\ref{eqsr1}) of the coherent projection $Q$ implies that: 
\begin{equation}
Q^{x}\,  \left(\exp\left( \sum_{i=1}^{n} a_{i} X_{i} \right)(x)\right) \  = \
\exp\left( \sum_{deg X_{i} = 1} a_{i}   X_{i} \right)(x)
\label{eqsr2}
\end{equation}
Therefore  $\displaystyle Q^{x}$ can be seen as a projection onto the 
(classical differential) geometric distribution. 
 
 \begin{remark}
The projection $\displaystyle Q^{x}$ has one more interesting feature: for any $x$ and 
$$\displaystyle u \, = \, \exp\left( \sum_{deg X_{i} = 1} a_{i} X_{i}
\right)(x)$$ 
we have $\displaystyle Q^{x} u \, = \, u$ and the curve 
$$s \in [0,1] \mapsto \delta^{x}_{s} \, u\ = \ \exp\left( s \, \sum_{deg X_{i} =
1} a_{i}   X_{i} \right)(x) $$
is $D$-horizontal and joins $x$ and $u$. This will be related to the
supplementary condition (B) further. 
\label{rkc}
\end{remark}

We may equally define a coherent projection which induces the dilatations 
$\delta$ from $\tilde{\delta}$.  Also, if we change the chosen normal frame with 
another of the same kind, we shall pass to a dilatation structure which is 
equivalent to $(U,d, \delta)$. In conclusion, coherent
projections are not geometrical objects per se, but in a natural way one may
define a notion of equivalent coherent projections such that the equivalence
class is geometrical, i.e. independent of the choice of a pair 
of particular dilatation structures, each in a given equivalence class. Another
way of putting this is that a class of equivalent dilatation structures may be
seen as a category and a coherent projection is a functor between such
categories. We shall not pursue this line here. 

The bottom line is that $\displaystyle (U, \bar{d}, \bar{\delta})$ is a
dilatation structure which belongs to an equivalence class which is independent
on the distribution $D$, and also independent on the choice of parameterization
$\phi$. It is associated to the manifold $M$ only. On the other hand  
$\displaystyle (U, \bar{d}, \bar{\delta})$ belongs to an equivalence class which
is depending only on the distribution $D$ and metric $g$ on $D$, thus intrinsic
to the sub-riemannian manifold $(M,D,g)$. The only advantage of choosing 
$\displaystyle \bar{\delta}, \delta$ related by the normal frame 
$\displaystyle \left\{ X_{1}, ..., X_{n} \right\}$ is that they are associated 
with a coherent projection with a simple expression.

\subsection{Length functionals associated to coherent projections}

\begin{definition}
Let $(X,\bar{d},\bar{\delta})$  be a strong dilatation structure with the Radon-Nikodym property  and 
$Q$ a coherent projection. We define the associated  distance $d: X \times X \rightarrow
[0, +\infty]$  by: 
$$d(x,y) \ = \ \inf \left\{ \int_{a}^{b} \bar{d}^{c(t)}(c(t), \dot{c}(t)) \mbox{ d}t \mbox{ :
} c: [a,b] \rightarrow X \mbox{ $\bar{d}$-Lipschitz },   \right.$$
$$\left. c(a) = x, c(b) = y, \mbox{ and }\forall a.e. \, \,  t \in [a,b] \quad  Q^{c(t)} \dot{c}(t) \ = \ \dot{c}(t)
\right\} $$
\end{definition}

The relation $x \equiv y$ if $d(x,y)< + \infty$ is an equivalence relation. The space
$X$ decomposes into a  reunion of equivalence classes, each equivalence class
being connected by horizontal curves. 

It is easy to see that $d$ is a finite  distance on each equivalence class. Indeed, 
from theorem \ref{fleng} we deduce that for any $x,y \in X$ $d(x,y) \geq \bar{d}(x,y)$.  
Therefore $d(x,y) = 0$ implies $x = y$. The other properties of a 
distance are straightforward. 

Later we shall give a sufficient condition (the generalized Chow condition (Cgen))
 on the coherent projection $Q$ for $X$ to be (locally) connected by horizontal curves.

\begin{proposition}
Suppose that $X$ is connected by
horizontal curves and $(X,d)$ is complete. Then  $d$ is  a length distance. 
\end{proposition}

\paragraph{Proof.}
Because $(X,d)$ is complete, it is sufficient 
to check that $d$ has the approximate middle property:  for any $\varepsilon > 0$ 
and for any $x, y \in X$ there exists $z \in X$ such that 
$$\displaystyle \max \left\{ d(x,z) , d(y,z) \right\} \leq \frac{1}{2} \, d(x,y) +
\varepsilon$$

Given $\varepsilon > 0$, from the definition of $d$ we deduce that there exists a
horizontal curve $c: [a,b] \rightarrow X$ such that $c(a) = x$, $c(b) = y$ and 
$\displaystyle d(x,y) + 2 \varepsilon \geq l(c)$ (where $l(c)$ is the length of $c$ 
with respect to the distance $\bar{d}$). There exists then $\tau \in [a,b]$ such that 
$$\int_{a}^{\tau} \bar{d}^{c(t)}(c(t), \dot{c}(t)) \mbox{ d}t \ = \ \int_{\tau}^{b} 
\bar{d}^{c(t)}(c(t), \dot{c}(t)) \mbox{ d}t \ = \ \frac{1}{2}\, l(c) $$
Let $z = c(\tau)$. We have then: 
$\displaystyle \max \left\{ d(x,z) , d(y,z) \right\} \ \leq \ \frac{1}{2} \, l(c) \ \leq \ 
\frac{1}{2} \, d(x,y) + \varepsilon$. 
Therefore $d$ is a length distance.
\hfill $\square$

{\bf Notations concerning length functionals.} The length functional associated to the distance $\bar{d}$ is denoted by 
$\bar{l}$. In the same way the length functional associated with 
$\displaystyle \bar{\delta}^{x}_{\varepsilon}$ is denoted by 
$\displaystyle \bar{l}^{x}_{\varepsilon}$. 
  
 We introduce the space $\displaystyle  \mathcal{L}_{\varepsilon}(X,d, \delta) \subset X 
 \times Lip([0,1],X,d)$: 
 $$\mathcal{L}_{\varepsilon}(X,d,  \delta) \ = \ \left\{(x ,c) \in X  \times
 \mathcal{C}([0,1],X) 
\mbox{ : } c: [0,1] \in U(x) \, \, , \,     \right.$$ 
$$\left. \delta^{x}_{\varepsilon}c \mbox{ is } 
\bar{d}-Lip, \quad Q-\mbox{horizontal} \mbox{ and } Lip(\delta^{x}_{\varepsilon}c) \leq 2 \varepsilon
l_{d}(\delta^{x}_{\varepsilon}c) 
\right\} $$

For any $\varepsilon \in (0,1)$ we define the length functional 
$$l_{\varepsilon}: \mathcal{L}_{\varepsilon}(X,d,  \delta) \rightarrow
[0,+\infty] \quad , \quad l_{\varepsilon}(x,c) \ = \ l^{x}_{\varepsilon}(c) \ = \ \frac{1}{\varepsilon}
\, 
\bar{l}(\delta^{x}_{\varepsilon} c)  $$

By theorem \ref{fleng} we have: 
$$l^{x}_{\varepsilon}(c) \ = \ \int^{1}_{0} \, \frac{1}{\varepsilon} \,
\bar{d}^{\delta^{x}_{\varepsilon} c(t)} \left( \delta^{x}_{\varepsilon} c(t), 
\frac{\bar{d}}{dt} \left(
\delta^{x}_{\varepsilon} c \right) (t)\right) \mbox{ d}t \ = $$
$$= \ \int^{1}_{0} \, \frac{1}{\varepsilon} \,
\bar{d}^{\delta^{x}_{\varepsilon} c(t)} \left( \delta^{x}_{\varepsilon} c(t), 
\delta^{x}_{\varepsilon} \, \frac{\hat{d}^{x}_{\varepsilon}}{dt} \, c(t)
\right) \mbox{ d}t $$
Another description of the length functional $\displaystyle l^{x}_{\varepsilon}$
is the following. 

\begin{proposition}
For any $\displaystyle (x,c) \in \mathcal{L}_{\varepsilon}(X,d,  \delta) $ we
have $$\displaystyle l^{x}_{\varepsilon}(c) \, = \, \bar{l}^{x}_{\varepsilon} 
\left(Q^{x}_{\varepsilon} c \right)  $$
\label{phor}
\end{proposition}

\paragraph{Proof.}
Indeed, we shall use an alternate definition of the length functional. Let 
$c$ be a curve such that $\displaystyle \delta^{x}_{\varepsilon} c$ is 
$\bar{d}$-Lipschitz and $Q$-horizontal. Then: 
$$l^{x}_{\varepsilon}(c) \ = \ \sup \left\{ \sum^{n}_{i=1} 
\frac{1}{\varepsilon} \bar{d}\left( \delta^{x}_{\varepsilon} c(t_{i}) , 
\delta^{x}_{\varepsilon} c(t_{i+1})\right) \mbox{ : } 0= t_{1} < ... < t_{n+1} =
1 \right\} \ = $$
$$= \ \sup \left\{ \sum^{n}_{i=1} 
\frac{1}{\varepsilon} \bar{d}\left( \bar{\delta}^{x}_{\varepsilon}
Q^{x}_{\varepsilon} c(t_{i}) , 
\bar{\delta}^{x}_{\varepsilon} Q^{x}_{\varepsilon} c(t_{i+1})\right) \mbox{ : } 0= t_{1} < ... < t_{n+1} =
1 \right\} \ = $$ 
$$= \ \bar{l}^{x}_{\varepsilon} 
\left(Q^{x}_{\varepsilon} c \right) $$
\hfill $\square$

\subsection{Supplementary hypotheses}

\begin{definition}
Let $\displaystyle (X, \bar{d}, \bar{\delta})$ be a strong dilatation structure 
and $Q$ a coherent projection. Further is a list of supplementary hypotheses on 
$Q$: 
\begin{enumerate}
\item[(A)] $\displaystyle \delta^{x}_{\varepsilon}$ is $\bar{d}$-bilipschitz 
in compact sets in the following sense: for any compact set $K \subset X$ and 
for any $\varepsilon \in (0,1]$ there is a number $L(K) > 0$ such that for any 
$x \in K$ and $u,v$ sufficiently close to $x$ we have: 
$$\frac{1}{\varepsilon} \, \bar{d} \left(\delta^{x}_{\varepsilon} u ,
\delta^{x}_{\varepsilon} v \right) \, \leq \, L(K) \, \bar{d}(u,v)$$  
\item[(B)] if $\displaystyle u = Q^{x} u$ then the curve $\displaystyle 
t \in [0,1] \mapsto \, Q^{x}\, \delta^{x}_{t} \, u \, = \, \bar{\delta}^{x}_{t}
u \, = \, \delta^{x}_{t} u$ is $Q$-everywhere horizontal and for any $a \in [0,1]$ we have 
$$  \limsup_{a \rightarrow 0} \frac{\bar{l}\left( t \in [0,a] \mapsto
 \bar{\delta}^{x}_{t} u \right)}{\bar{d}(x, \bar{\delta}^{x}_{a} u)} \, =   \, 1$$ 
 uniformly with respect to  $x,u$ in compact set $K$. 
 \end{enumerate}
 \end{definition}
 
 Condition  (A), as well as the  property (IV) definition  
\ref{defcoh}, is another smoothness condition  on $Q$ with
 respect to the strong dilatation structure $\displaystyle (X , \bar{d},
 \bar{\delta})$. 
 
 The condition (A)  has several useful consequences, among them the fact that 
 for any $\bar{d}$-Lipschitz curve $c$, the curve $\displaystyle
 \delta^{x}_{\varepsilon}c$ is also Lipschitz. Another consequence is that 
 $\displaystyle Q^{x}_{\varepsilon}$ is locally $\bar{d}$-Lipschitz.  
 More precisely,  for any compact set $K \subset X$ and 
for any $\varepsilon \in (0,1]$ there is a number $L(K) > 0$ such that for any 
$x \in K$ and $u,v$ sufficiently close to $x$ we have:
\begin{equation}
\left(\bar{\delta}^{x}_{\varepsilon} \bar{d}\right) \, 
\left(Q^{x}_{\varepsilon} u , Q^{x}_{\varepsilon} v \right) \, \leq \, L(K) \, \bar{d}(u,v)
\label{uu}
\end{equation}
with the notation 
$$\left(\bar{\delta}^{x}_{\varepsilon} \bar{d} \right) (u, v) \, = \, 
\frac{1}{\varepsilon} \, \bar{d} \left( \bar{\delta}^{x}_{\varepsilon} u , 
\bar{\delta}^{x}_{\varepsilon} v \right)$$
Indeed, we have: 
$$\left(\bar{\delta}^{x}_{\varepsilon} \bar{d}\right) 
\left(Q^{x}_{\varepsilon} u , Q^{x}_{\varepsilon} v \right) \, = \, 
\frac{1}{\varepsilon} \, \bar{d} \left(\delta^{x}_{\varepsilon} u ,
\delta^{x}_{\varepsilon} v \right) \, \leq \, L(K) \, \bar{d}(u,v)$$

See the remark \ref{rkc} for the meaning of the  condition B for the case
sub-riemannian geometry, where it is explained why condition B is a 
generalization of the fact that   the
"distribution" $\displaystyle x \mapsto Q^{x} U(x)$ is  generated by horizontal  one 
parameter flows.

 Condition (B) will be useful later, along with the generalized Chow 
  condition (Cgen).

\section{The generalized Chow condition}
\label{sechormander}

{\bf Notations about words.} For any set $A$ we denote by $\displaystyle A^{*}$ the collection of 
finite words $\displaystyle q = a_{1}...a_{p}$, $\displaystyle p \in 
\mathbb{N}$, $p > 0$. The empty word is denoted 
by $\emptyset$. The length of the word $\displaystyle q = a_{1}...a_{p}$ 
is $\mid q \mid = p$; the length of the empty word is $0$. 

The collection of words infinite at right over the alphabet $A$ is denoted by 
$\displaystyle A^{\omega}$. For any word $\displaystyle w \in A^{\omega} \cup
A^{*}$ and
any $p \in \mathbb{N}$ we denote by $\displaystyle [w]_{p}$ the finite word 
obtained from the first $p$ letters of $w$ (if $p=0$ then $\displaystyle 
[w]_{0} = \emptyset$ (in the case of a finite word $q$, if $p > \mid q \mid$ 
then $\displaystyle [q]_{p} = q$).

For any non-empty $\displaystyle q_{1}, q_{2} \in A^{*}$ and $w \in A^{\omega}$ the 
concatenation  of $\displaystyle q_{1}$ and $\displaystyle q_{2}$ is the 
finite word $\displaystyle q_{1}q_{2} \in A^{*}$ and the concatenation of 
$\displaystyle q_{1}$ and $w$ is the (infinite) word $\displaystyle q_{1}w \in 
A^{\omega}$. The empty word $\emptyset$ is seen both as an infinite 
word or a finite word and for any  $\displaystyle q \in A^{*}$ and 
$\displaystyle w \in A^{\omega}$ we have 
$\displaystyle q \emptyset = q$ (as concatenation of finite words) and 
$\emptyset w = w$ (as concatenation of a finite empty word and an infinite 
word).

\subsection{Coherent projections as transformations of words} To any coherent projection $Q$ in a strong dilatation structure $\displaystyle 
(X, \bar{d}, \bar{\delta})$ we associate a family of transformations as 
follows. 

\begin{definition}
For any non-empty word $\displaystyle w \in (0,1]^{\omega}$ and any 
$\varepsilon \in (0,1]$ we define the transformation 
$$\Psi_{\varepsilon w} : X^{*}_{\varepsilon w} \, \subset X^{*} \setminus
\left\{ \emptyset \right\} \, 
\rightarrow X^{*} $$  
given by: for any non-empty finite word $\displaystyle q = x x_{1} ... x_{p} \in 
X^{*}_{\varepsilon w}$ we have 
$$ \Psi_{\varepsilon w} (x x_{1} ... x_{p}) \ = \ 
\Psi_{\varepsilon w}^{1}(x) ... \Psi_{\varepsilon w}^{k+1}(x x_{1} ... x_{k}) ... 
\Psi_{\varepsilon w}^{p+1}(x x_{1} ... x_{p}) $$  
The functions $\displaystyle \Psi_{\varepsilon w}^{k}$ are defined by: 
 $\displaystyle \Psi_{\varepsilon w}^{1}(x) = x$, and for any $k \geq 1$ we 
 have 
 \begin{equation}
\Psi_{\varepsilon w}^{k+1}([q]_{k+1}) \ = \
\delta^{x}_{\varepsilon^{-1}} \, Q_{w_{k}}^{\delta^{x}_{\varepsilon} \,
\Psi_{\varepsilon w}^{k}([q]_{k})} \, \delta_{\varepsilon}^{x} \, q_{k+1} 
\label{recrel}
\end{equation}
If $w = \emptyset$ then $\displaystyle \Psi_{\varepsilon \emptyset}^{k}$ is
defined as previously $\displaystyle \Psi_{\varepsilon \emptyset}^{1}(x) = x$, 
with the only difference that for any $k \geq 1$ we 
 have 
$$\Psi_{\varepsilon \emptyset}^{k+1}([q]_{k+1}) \ = \
\delta^{x}_{\varepsilon^{-1}} \, Q^{\delta^{x}_{\varepsilon} \,
\Psi_{\varepsilon w}^{k}([q]_{k})} \, \delta_{\varepsilon}^{x} \, q_{k+1} $$
\label{dwords}
\end{definition}

The domain $\displaystyle X^{*}_{\varepsilon w} \subset X^{*} \setminus
\left\{ \emptyset \right\}$ is such that the previous definition makes sense. By
using the definition of a coherent projection, we may redefine $\displaystyle
X^{*}_{\varepsilon w}$ as  follows:  for any compact set $K \subset X$  there 
is $\rho = \rho(K) > 0$ such that for any $x \in K$ the word 
$\displaystyle q = x x_{1} ... x_{p} \, \in X^{*}_{\varepsilon w}$ if for any 
$k \geq 1$ we have 
$$\bar{d}\left( x_{k+1} , \Psi_{\varepsilon w}^{k}([q]_{k}) \right) \, \leq \,
\rho $$

We shall explain the meaning of these transformations for 
$\varepsilon = 1$. 

\begin{proposition}
Suppose that condition (B) holds for the coherent projection $Q$. If 
$$\displaystyle y \ = \ \Psi_{1\emptyset}^{k+1}(x x_{1} ... x_{k})$$ 
then there is a $Q$-horizontal curve joining $x$ and $y$.
\end{proposition}

\paragraph{Proof.} 
 By  definition \ref{dwords} for 
$\varepsilon = 1$  we have: 
$$\Psi_{1w}^{1}(x) =  x \quad , \quad \Psi_{1w}^{2}(x, x_{1}) \ =  \
Q_{w_{1}}^{x} \, x_{1} \quad , \quad$$ 
$$\Psi_{1w}^{3}(x, x_{1}, x_{2}) \ = \ 
Q_{w_{2}}^{Q_{w_{1}}^{x} x_{1}} \, x_{2} \quad ... $$
Suppose now that condition (B) holds for the coherent projection $Q$. Then the
curve $\displaystyle t \in [0,1] \mapsto \bar{\delta}^{x}_{t} Q^{x} u$ is a 
$Q$-horizontal curve joining $x$ with $\displaystyle Q^{x} u$. Therefore 
 by applying inductively the condition (B) we get that 
 there is a $Q$-horizontal curve between 
$\displaystyle \Psi_{1\emptyset}^{k}(x x_{1} ... x_{k-1})$ and $\displaystyle
\Psi_{1\emptyset}^{k+1}(x x_{1} ... x_{k})$ for any $k > 1$ and a $Q$-horizontal
curve joining $x$ and $\displaystyle \Psi_{1\emptyset}^{2}(x x_{1})$. 
\hfill $\square$

There are three more properties of the transformations $\displaystyle
\Psi_{\varepsilon w}$. 

\begin{proposition}
With the notations from definition \ref{dwords} we have: 
\begin{enumerate}
\item[(a)] $\displaystyle \Psi_{\varepsilon w} \, 
\Psi_{\varepsilon \emptyset} \, = \, \Psi_{\varepsilon \emptyset}$. Therefore 
we have the equality of sets: 
$$\Psi_{\varepsilon \emptyset}\, \left( X^{*}_{\varepsilon \emptyset} \cap x X^{*}
\right) \, = \, \Psi_{\varepsilon w}\, \left( \Psi_{\varepsilon \emptyset}\,
\left( X^{*}_{\varepsilon \emptyset} \cap x X^{*}
\right) \right)$$
\item[(b)] $\displaystyle \Psi_{\varepsilon \emptyset}^{k+1} (x q_{1} ... q_{k})
\, = \, \delta_{\varepsilon^{-1}}^{x} \, \Psi_{1 \emptyset}^{k+1} (x
\delta_{\varepsilon}^{x}q_{1} ... \delta_{\varepsilon}^{x} q_{k})$
\item[(c)] $\displaystyle \lim_{\varepsilon \rightarrow 0} 
 \delta_{\varepsilon^{-1}}^{x} \, \Psi_{1 \emptyset}^{k+1} (x
\delta_{\varepsilon}^{x}q_{1} ... \delta_{\varepsilon}^{x} q_{k}) \, = \, 
\Psi_{0 \emptyset}^{k+1} (x q_{1} ... q_{k})$ uniformly with respect to 
$\displaystyle x, q_{1}, ..., q_{k}$ in compact set. 
\end{enumerate}
\end{proposition}

\paragraph{Proof.}
(a) We use induction on $k$ to prove that for any natural number $k$ we have: 
\begin{equation}
\Psi^{k+1}_{\varepsilon w} \left( \Psi^{1}_{\varepsilon \emptyset}(x) ... 
\Psi^{k+1}_{\varepsilon \emptyset} (x q_{1} ... q_{k}) \right) \, = \, 
\Psi^{k+1}_{\varepsilon \emptyset} (x q_{1} ... q_{k})
\label{prpaproof}
\end{equation}
 For $k = 0$ we have have to prove that  $x = x$ which is trivial. For 
 $k = 1$ we have to prove that 
 $$\Psi^{2}_{\varepsilon w} \left( \Psi^{1}_{\varepsilon \emptyset}(x) \,  
\Psi^{2}_{\varepsilon \emptyset} (x q_{1}) \right) \, = \, 
\Psi^{2}_{\varepsilon \emptyset} (x q_{1})$$ 
This means: 
$$\Psi^{2}_{\varepsilon w} \left( x \,  \delta_{\varepsilon^{-1}}^{x} \, 
Q^{x} \delta^{x}_{\varepsilon} q_{1} \right) \, = \, 
\delta_{\varepsilon^{-1}}^{x} \, Q_{w_{1}}^{x} \, \delta_{\varepsilon}^{x} \, 
\delta^{x}_{\varepsilon^{-1}} \, Q^{x} \, \delta^{x}_{\varepsilon} x_{1} \, = 
$$ 
$$= \, \delta^{x}_{\varepsilon^{-1}} \, Q^{x} \, \delta^{x}_{\varepsilon} x_{1} \, =
\, \Psi^{2}_{\varepsilon \emptyset} (x q_{1})$$

  Suppose now that $l \geq 2$ and  for any $k \leq l$ the relations 
(\ref{prpaproof}) are true. Then, as previously, it is easy to check 
(\ref{prpaproof}) for $k = l+1$.

(b) is true by direct computation. The point (c) is a straightforward
consequence of (b) and definition of coherent projections.  
\hfill $\square$

\begin{definition}
Let $N \in \mathbb{N}$ be a strictly positive natural number and $\varepsilon
\in (0,1]$. 
We say that $x \in X$ is {\bf $(\varepsilon,N,Q)$-nested} in a open neighbourhood 
$U \subset X$ if there is $\rho> 0$ such that for any finite word $\displaystyle 
q = x_{1} ... x_{N} \in X^{N}$ with 
$$\bar{\delta}^{x}_{\varepsilon} \bar{d} \left( x_{k+1} , \Psi_{\varepsilon
 \emptyset}^{k}([xq]_{k}) \right) \, \leq \, \rho $$
for any $k = 1, ... , N$, we have $\displaystyle q \in U^{N}$. 

If $x \in U$ is $(\varepsilon,N,Q)$-nested then denote by 
$\displaystyle U(x,\varepsilon, N,Q,\rho) \subset U^{N}$ the collection of 
words $\displaystyle q \in U^{N}$ such that 
$\displaystyle \bar{\delta}^{x}_{\varepsilon} \bar{d}\left( x_{k+1} ,
\Psi_{\varepsilon \emptyset}^{k}([xq]_{k}) \right)
\, < \,
\rho $ for any $k = 1, ... , N$. 
\end{definition}

\begin{definition}
A coherent projection $Q$  satisfies  the {\bf generalized Chow condition} 
if: 
\begin{enumerate}
 \item[(Cgen)] for any compact set $K$ there are $\rho = \rho(K) > 0$, $r = r(K) > 0$, 
 a natural number $N = N(Q,K)$ and a function 
 $F(\eta) = \mathcal O(\eta)$ such that for
 any $x \in K$ and $\varepsilon \in (0,1]$ 
 there are   neighbourhoods $U(x)$, $V(x)$ such that  any 
  $x \in K$ is $(\varepsilon, N,Q)$-nested in $U(x)$, 
  $\displaystyle B(x,r, \bar{\delta}^{x}_{\varepsilon} \bar{d}) 
  \subset V(x)$  and such that the mapping 
 $$x_{1} ... x_{N} \in U(x,N,Q,\rho) \, \mapsto \,  
 \Psi^{N+1}_{\varepsilon \emptyset}(x x_{1} ... x_{N})$$ 
 is surjective from $\displaystyle U(x,\varepsilon,N,Q,\rho)$ to 
 $\displaystyle V(x)$. 
 Moreover for any $\displaystyle z \in V(x)$ there exist $\displaystyle  
 y_{1}, ...  y_{N} \in U(x,\varepsilon, N, Q, \rho)$ such that 
 $\displaystyle z \, = \, \Psi^{N+1}_{\varepsilon \emptyset}(x y_{1}, ...  y_{N})$ and 
 for any $k = 0, ... , N-1$ we have 
 $$\delta^{x}_{\varepsilon} \bar{d} \left( \Psi^{k+1}_{\varepsilon
 \emptyset}(x y_{1} ...  y_{k}) , 
 \Psi^{k+2}_{\varepsilon \emptyset}(x y_{1} ...  y_{k+1}) \right) \ \leq \
 F(\delta^{x}_{\varepsilon} \bar{d}(x,z))$$
 \end{enumerate}
 \label{defhorgen} 
 \end{definition}

Condition  (Cgen) is inspired from lemma 1.40 Folland-Stein \cite{folstein}. If 
the coherent projection $Q$ satisfies also (A) and (B) then   in the
space $\displaystyle (U(x), \bar{\delta}^{x}_{\varepsilon})$, with 
coherent projection $\displaystyle \hat{Q}^{x, \cdot}_{\varepsilon.\cdot}$, 
 we can join any two sufficiently close points by a sequence of 
at most $N$ horizontal curves. Moreover
there is a control on the length of these curves via condition (B) and condition
(Cgen); in sub-riemannian geometry the function $F$ is of the type 
$\displaystyle   F(\eta) = \eta^{1/m}$ with $m$ positive natural
number.

\begin{definition}
Suppose that the coherent projection $Q$ satisfies conditions (A), (B) and (Cgen).  
Let us consider $\varepsilon \in (0,1]$ and 
 $x, y \in K$, $K$ compact in $X$. With the notations from 
  definition \ref{defhorgen},  
  suppose that there are numbers $N=N(Q,K)$, $\rho = \rho(Q,K) > 0$  and 
  words $\displaystyle  x_{1} ... x_{N} \in U(x,\varepsilon,N,Q,\rho)$
  such that  
  $$ y = \Psi_{\varepsilon\emptyset}^{N+1}(x x_{1} ... x_{N})  $$
  To these data we associate a {\bf short curve}  joining 
  $x$ and $y$,  $c : [0,N] \rightarrow X$ defined by: 
 for any  $t \in [0, N]$ then let $k = [t]$, where $[b]$ is the integer
part of the real number $b$. We define the short curve by 
$$c(t) \ = \ \bar{\delta}^{x, \Psi_{\varepsilon \emptyset}^{k+1}(x x_{1} ...
x_{k})}_{\varepsilon,t+N -k} 
Q^{\Psi_{\varepsilon\emptyset}^{k+1}(x   x_{1} ... x_{k})} x_{k+1} $$
Any short  curve joining $x$ and $y$ is a increasing linear reparameterization of 
a curve $c$ described previously.
\end{definition}


\subsection{The candidate tangent space}
\label{candidate}

Let $(X, \bar{d}, \bar{\delta})$ be a strong dilatation structure and 
$Q$   a coherent projection.  Then we have
the induced dilatations 
$$\mathring{\delta}^{x,u}_{\mu} v \ = \ \Sigma^{x}(u, \delta^{x}_{\mu}
\Delta^{x}(u, v))$$
and the induced  projection 
$$\mathring{Q}^{x,u}_{\mu} v \ = \ \Sigma^{x}(u, Q^{x}_{\mu}
\Delta^{x}(u, v)) $$
For any curve $c: [0,1] \rightarrow U(x)$ which is 
$\displaystyle \mathring{\delta}^{x}$-derivable and 
$\displaystyle \mathring{Q}^{x}$-horizontal almost everywhere: 
$$\frac{\mathring{d}^{x}}{dt} c(t) \ = \ \mathring{Q}^{x,u} \, 
\frac{\mathring{d}^{x}}{dt} c(t)$$
we define the length
$$l^{x}(c) \ = \ \int^{1}_{0} \, \bar{d}^{x}\left( x, \Delta^{x}(c(t), 
\frac{\mathring{d}^{x}}{dt} c(t)) \right) \mbox{ d}t $$
and the distance function: 
$$\mathring{d}^{x}(u,v) \, = \, \inf \left\{ l^{x}(c) \mbox{ : }  c: [0,1]
\rightarrow U(x) \, \mbox{ is $\displaystyle \mathring{\delta}^{x}$-derivable} 
,\right. $$ 
$$\left. \mbox{ and $\displaystyle \mathring{Q}^{x}$-horizontal a.e.} \, , 
\, c(0) = u, c(1) = v \right\}$$

We want to prove that  $\displaystyle (U(x), \mathring{d}^{x}, \mathring{\delta}^{x})$ is a
strong dilatation structure and $\displaystyle \mathring{Q}^{x}$ is a coherent projection. 
For this we need first the following proposition.

\begin{proposition}
The curve $c: [0,1] \rightarrow U(x)$ is 
$\displaystyle \mathring{\delta}^{x}$-derivable, 
$\displaystyle \mathring{Q}^{x}$-horizontal almost everywhere, and 
$\displaystyle l^{x}(c) < + \infty$ if and only if the curve  
$\displaystyle Q^{x} c$ is $\displaystyle \bar{\delta}^{x}$-derivable 
almost everywhere and 
$\displaystyle \bar{l}^{x}(Q^{x} c) < + \infty$. Moreover, we have  
$$\displaystyle \bar{l}^{x}(Q^{x} c) \, = \, l^{x}(c)$$ 
\label{propcan}
\end{proposition}

\paragraph{Proof.}
The curve $c$ is $\displaystyle \mathring{Q}^{x}$-horizontal almost everywhere
if and only if for almost any $t \in [0,1]$ we have 
$$Q^{x} \, \Delta^{x}(c(t) , \frac{\mathring{d}^{x}}{dt} c(t)) \ = \ 
\Delta^{x}(c(t) , \frac{\mathring{d}^{x}}{dt} c(t)) $$
We shall prove that  
 $c$ is $\displaystyle \mathring{Q}^{x}$-horizontal  is
equivalent with 
\begin{equation}
\Theta^{x}(c(t), 
\frac{\mathring{d}^{x}}{dt} c(t)) \ = \  \frac{\bar{d}^{x}}{dt} \left(Q^{x} c
\right)(t) 
\label{blakdot}
\end{equation}
Indeed, \eqref{blakdot} is equivalent with
$$\lim_{\varepsilon \rightarrow 0} \bar{\delta}^{x}_{\varepsilon^{-1}} 
\bar{\Delta}^{x}(Q^{x} c(t), Q^{x} c(t+\varepsilon)) \ = \ 
\bar{\Delta}^{x}(Q^{x} c(t) , \Theta^{x}(c(t), 
\frac{\mathring{d}^{x}}{dt} c(t)))$$ 
which is equivalent with 
$$\lim_{\varepsilon \rightarrow 0} \bar{\delta}^{x}_{\varepsilon^{-1}} 
\bar{\Delta}^{x}(Q^{x} c(t), Q^{x} c(t+\varepsilon)) \ = \ 
\Delta^{x}(c(t), \frac{\mathring{d}^{x}}{dt} c(t))$$ 
But this is equivalent with: 
\begin{equation}
\lim_{\varepsilon \rightarrow 0} \bar{\delta}^{x}_{\varepsilon^{-1}} 
\bar{\Delta}^{x}( Q^{x} c(t), Q^{x} c(t+\varepsilon)) \ = \ 
\lim_{\varepsilon \rightarrow 0} \delta^{x}_{\varepsilon^{-1}} 
\Delta^{x}(c(t), c(t+ \varepsilon))
\label{blakdot1}
\end{equation}
The horizontality condition for the curve $c$ can be written as: 
$$\lim_{\varepsilon \rightarrow 0} Q^{x} \delta^{x}_{\varepsilon^{-1}} 
\Delta^{x}(c(t), c(t+\varepsilon)) \ = \ \lim_{\varepsilon \rightarrow 0}  
\delta^{x}_{\varepsilon^{-1}} 
\Delta^{x}(c(t), c(t+\varepsilon))$$ 
We use now the properties of $\displaystyle Q^{x}$ in the left hand side of the
previous equality: 
$$Q^{x} \delta^{x}_{\varepsilon^{-1}} 
\Delta^{x}(c(t), c(t+\varepsilon)) \ = \ 
\bar{\delta}^{x}_{\varepsilon^{-1}} Q^{x} \Delta^{x}(c(t), c(t+\varepsilon)) 
\ =$$
$$= \  \bar{\delta}^{x}_{\varepsilon^{-1}} \bar{\Delta}^{x}(Q^{x} c(t), Q^{x}
c(t+\varepsilon))$$
thus after taking the limit as $\varepsilon \rightarrow 0$ we prove that 
the limit $$\displaystyle \lim_{\varepsilon \rightarrow 0} 
\bar{\delta}^{x}_{\varepsilon^{-1}} \bar{\Delta}^{x}(Q^{x} c(t), Q^{x}
c(t+\varepsilon))$$ exists and we obtain: 
$$\lim_{\varepsilon \rightarrow 0}  
\delta^{x}_{\varepsilon^{-1}} 
\Delta^{x}(c(t), c(t+\varepsilon)) \ = \  \lim_{\varepsilon \rightarrow 0} 
\bar{\delta}^{x}_{\varepsilon^{-1}} \bar{\Delta}^{x}(Q^{x} c(t), Q^{x}
c(t+\varepsilon))$$
This last equality is the same as \eqref{blakdot1}, which is equivalent with 
\eqref{blakdot}. 

As a consequence we obtain the following equality, for almost any $t \in [0,1]$:
\begin{equation}
\bar{d}^{x}\left( x, \Delta^{x}(c(t), 
\frac{\mathring{d}^{x}}{dt} c(t)) \right) \ = \ 
\bar{\Delta}^{x}(Q^{x} c(t) , \frac{\bar{d}^{x}}{dt} \left(Q^{x} c
\right)(t)) 
\label{blakdot3}
\end{equation}
This implies that $Q^{x}c$ is absolutely continuous and by theorem 
\ref{tupper}, as in the proof of theorem \ref{fleng} (but without using 
the Radon-Nikodym property  property, because we already know that $\displaystyle Q^{x}c$ is
derivable a.e.), we obtain the following formula for the length of the curve 
$Q^{x}c$: 
$$\bar{l}^{x}(Q^{x} c) \ = \ \int^{1}_{0} \, \bar{d}^{x}\left( 
x, , \bar{\Delta}^{x}(Q^{x} c(t) , \frac{\bar{d}^{x}}{dt} \left(Q^{x} c
\right)(t)) \right) \mbox{ d}t  $$
But we have also:
$$l^{x}(c) \ = \ \int^{1}_{0} \, \bar{d}^{x}\left( x, \Delta^{x}(c(t), 
\frac{\mathring{d}^{x}}{dt} c(t)) \right) \mbox{ d}t  $$
By \eqref{blakdot3} we obtain $\displaystyle \bar{l}^{x}(Q^{x} c) \, = \,
l^{x}(c)$. 
\hfill $\square$

\begin{proposition}
If $(X, \bar{d}, \bar{\delta})$ is a strong dilatation structure, 
$Q$ is  a coherent projection  and $\mathring{d}^{x}$ is finite then 
 the triple  $\displaystyle (U(x), \Sigma^{x}, \delta^{x})$
is a normed conical group, with the norm induced by the  left-invariant
distance $\mathring{d}^{x}$.
\label{prevprop}
\end{proposition}

\paragraph{Proof.}
The fact that  $\displaystyle (U(x), \Sigma^{x}, \delta^{x})$ is a conical group
comes directly from the definition \ref{defcoh} of a coherent projection.
Indeed,  it is enough to use proposition \ref{p1proj} (c) and the formalism 
of binary decorated trees in \cite{buligadil1} section 4 (or theorem 11
\cite{buligadil1}), in order to  reproduce 
the part of the proof of theorem 10 (p.87-88) in that paper, concerning the
conical group structure. There is one small subtlety though. 
In the proof of theorem \ref{tgene}(a) the same modification of proof has
been done starting from the axiom A4+, namely the existence of the uniform limit
$\displaystyle \lim_{\varepsilon \rightarrow 0} \Sigma^{x}_{\varepsilon}(u,v) \,
= \, \Sigma^{x}(u,v)$. Here we need first to prove this limit, in a similar way
as in the corollary 9 \cite{buligadil1}. We shall use for this the distance 
$\displaystyle \mathring{d}^{x}$ instead of the distance 
in the metric tangent space of $(X,d)$ at $x$ denoted by 
 $\displaystyle d^{x}$ (which is not yet proven to exist). The distance 
 $\displaystyle \mathring{d}^{x}$ is supposed to be finite by hypothesis.
 Moreover, by its definition and proposition \ref{propcan} we have 
 $$\mathring{d}^{x}(u,v) \, \geq \, \bar{d}^{x}(u,v)$$
 therefore the distance $\displaystyle \mathring{d}^{x}$ is non degenerate. By 
 construction  this distance is also left invariant with respect to the 
 group operation $\displaystyle \Sigma^{x}(\cdot, \cdot)$. Therefore we may
 repeat the proof of corollary 9 \cite{buligadil1} and obtain the result that 
 A4+ is true for $(X,d,\delta)$. 

What we need
to prove next  is that $\mathring{d}^{x}$ induces a norm on the conical group 
 $\displaystyle (U(x), \Sigma^{x}, \delta^{x})$. For this it is enough to 
 prove that 
 \begin{equation}
 \mathring{d}^{x}(\mathring{\delta}^{x,u}_{\mu} v, 
 \mathring{\delta}^{x,u}_{\mu} w ) \, = \, \mu \, \mathring{d}^{x}(v,w)
 \label{mewn1}
 \end{equation}
 for any $v, w \in U(x)$. This is a direct consequence of relation 
 (\ref{blakdot3}) from the proof of the  proposition 
 \ref{propcan}. Indeed, by direct computation we get that for any curve 
 $c$ which is $\mathring{Q}^{x}$-horizontal a.e. we have: 
 $$l^{x}(\mathring{\delta}^{x,u}_{\mu} c) \, = \, 
 \int_{0}^{1} \bar{d}^{x} \left( x, \Delta^{x}\left(\mathring{\delta}^{x,u}_{\mu} c(t) , 
 \frac{\mathring{d}^{x}}{dt} \left( \mathring{\delta}^{x,u}_{\mu} c \right) (t)
 \right)\right) \mbox{ dt} \, = \, $$
$$= \, \int_{0}^{1} \bar{d}^{x} \left( x, \delta^{x}_{\mu} \Delta^{x}\left( c(t) , 
 \frac{\mathring{d}^{x}}{dt}  c (t)
 \right)\right) \mbox{ dt}$$
 But  $c$ is $\mathring{Q}^{x}$-horizontal a.e., which implies, via
 (\ref{blakdot3}),  that 
 $$\delta^{x}_{\mu} \Delta^{x}\left( c(t) , 
 \frac{\mathring{d}^{x}}{dt}  c (t) \right) \, = \, \bar{\delta}^{x}_{\mu} \Delta^{x}\left( c(t) , 
 \frac{\mathring{d}^{x}}{dt}  c (t)
 \right)$$ 
 therefore we have 
 $$l^{x}(\mathring{\delta}^{x,u}_{\mu} c) \, = \, \int_{0}^{1} \bar{d}^{x} \left( x,
 \bar{\delta}^{x}_{\mu} \Delta^{x}\left( c(t) , 
 \frac{\mathring{d}^{x}}{dt}  c (t)
 \right)\right) \mbox{ dt} \, = \, \mu \, l^{x}(c)$$ 
 This implies (\ref{mewn1}), therefore the  proof is done.
\hfill $\square$

\begin{theorem}
If the generalized Chow condition (Cgen) and condition (B) are true then  
$\displaystyle (U(x), \Sigma^{x}, \delta^{x})$ is local conical group which is 
a neighbourhood of the neutral element of a  Carnot group 
generated by $\displaystyle Q^{x} U(x)$. 
\label{tancarnot}
\end{theorem}

\paragraph{Proof.}
 For any 
$\varepsilon \in (0,1]$, as a consequence of proposition  \ref{pdesc} we can put 
the recurrence relations \eqref{recrel} in the form:
 \begin{equation}
\Psi_{\varepsilon w}^{k+1}([q]_{k+1}) \ = \ \Sigma^{x}_{\varepsilon} \left( 
\Psi_{\varepsilon w}^{k}([q]_{k}), Q_{w_{k}}^{\delta^{x}_{\varepsilon} \,
\Psi_{\varepsilon w}^{k}([q]_{k})} \, \Delta^{x}_{\varepsilon} \left( 
\Psi_{\varepsilon w}^{k}([q]_{k}), q_{k+1}\right) \right) 
\label{recrelnew}
\end{equation}
This recurrence relation allows us to prove by induction that for any 
$k$ the limit 
$$\Psi_{ w}^{k}([q]_{k}) \ = \ \lim_{\varepsilon \rightarrow 0} 
\Psi_{\varepsilon w}^{k}([q]_{k})$$ 
exists and it satisfies the recurrence relation: 
 \begin{equation}
\Psi_{0 w}^{k+1}([q]_{k+1}) \ = \ \Sigma^{x} \left( 
\Psi_{0 w}^{k}([q]_{k}), Q_{w_{k}}^{x} \, \Delta^{x} \left( 
\Psi_{0 w}^{k}([q]_{k}), q_{k+1}\right) \right) 
\label{recrelnew2}
\end{equation}
and the initial condition $\displaystyle \Psi_{0 w}^{1}(x) = x$. 
We pass to the limit in the generalized Chow condition (Cgen) and we thus 
obtain that a neighbourhood of the neutral element $x$ is (algebraically) 
generated by $\displaystyle Q^{x} U(x)$. Then the distance 
$\displaystyle \mathring{d}^{x}$. Therefore by proposition \ref{prevprop}
$\displaystyle (U(x), \Sigma^{x}, \delta^{x})$ is a normed conical group 
generated by $\displaystyle Q^{x} U(x)$. 

Let $c:[0,1] \rightarrow U(x)$ be the curve $\displaystyle c(t) = \delta^{x}_{t}
u$, with $\displaystyle u \in Q^{x} U(x)$. Then we have $\displaystyle Q^{x} c(t) = c(t) = \bar{\delta}^{x}_{t} u$.
From condition (B) we get that $c$ is $\bar{\delta}$-derivable at $t=0$. A short
computation of this derivative shows that: 
$$\frac{d \bar{\delta}}{dt} c(0) \, = \, u$$
Another easy computation shows that the curve $c$ is 
$\displaystyle \bar{\delta}^{x}$-derivable if and only if the curve $c$ is 
$\bar{\delta}$-derivable at $t=0$, which is true, therefore $c$ is 
$\displaystyle \bar{\delta}^{x}$-derivable, in particular at $t=0$. 
Moreover, the expression of the 
$\displaystyle \bar{\delta}^{x}$-derivative of $c$ shows that $c$ is also 
$\displaystyle Q^{x}$-everywhere horizontal (compare with the remark \ref{rkc}).
We use the proposition \ref{propcan} and  relation 
(\ref{blakdot}) from its
proof to deduce that $\displaystyle c = Q^{x} c$  is
$\mathring{\delta}^{x}$-derivable at $t=0$, thus for any 
$\displaystyle u \in Q^{x} U(x)$  and small enough $t, \tau \in (0,1)$ 
we have 
\begin{equation}
\mathring{\delta}^{x, x}_{t + \tau} u \, = \,
\bar{\Sigma}^{x}(\bar{\delta}^{x}_{t} u , \bar{\delta}^{x}_{\tau} u)
\label{need22}
\end{equation}
By previous proposition \ref{prevprop} and corollary 6.3 \cite{buligadil2} (here
proposition \ref{cor63}) 
the normed conical group $\displaystyle (U(x), \Sigma^{x}, \delta^{x})$ is in
fact locally a homogeneous group, i.e. a 
simply connected Lie group which admits a positive graduation given by the 
eigenspaces of $\displaystyle \delta^{x}$. Indeed, corollary 6.3
\cite{buligadil1} is originally about strong dilatation structures, but the
generalized Chow condition implies that the distances $d$, $\bar{d}$ and 
$\displaystyle \mathring{d}^{x}$ induce the same uniformity, which, along with 
proposition \ref{prevprop}, are the only 
things needed for the proof of this corollary.  The conclusion of corollary 
6.3 \cite{buligadil2} therefore is true, that is 
$\displaystyle (U(x), \Sigma^{x}, \delta^{x})$ is  locally a homogeneous 
group. Moreover it is 
 locally Carnot if and only if on the generating space $\displaystyle Q^{x} U(x)$ any 
dilatation $\displaystyle \mathring{\delta}^{x, x}_{\varepsilon} u \, = \, 
\bar{\delta}^{x}_{\varepsilon}$ is linear in $\varepsilon$. But this is true, as 
shown by relation (\ref{need22}). This ends the proof. 
\hfill $\square$

\subsection{Coherent projections induce length dilatation structures}
\label{subscls}

\begin{theorem}
If $\displaystyle (X,\bar{d}, \bar{\delta})$ is a  tempered strong dilatation 
structure, has the Radon-Nikodym property  and $Q$ is a coherent projection, which  
 satisfies (A), (B), (Cgen) then $\displaystyle (X, d, \delta)$ is a length 
 dilatation structure. 
\label{mainsrthm}
\end{theorem}

\paragraph{Proof.}

 We shall prove that: 
\begin{enumerate}
\item[(a)] for any function 
$\displaystyle \varepsilon \in (0,1) \mapsto (x_{\varepsilon}, c_{\varepsilon}) 
\in \mathcal{L}_{\varepsilon}(X,d,  \delta)$ which
converges to $\displaystyle (x,c)$ as $\varepsilon \rightarrow 0$, with $c: [0,1] \rightarrow U(x)$ 
$\displaystyle \mathring{\delta}^{x}$-derivable and 
$\displaystyle \mathring{Q}^{x}$-horizontal almost everywhere, we have: 
$$l^{x}(c) \ \leq \ \liminf_{\varepsilon \rightarrow 0}
l^{x_{\varepsilon}}(c_{\varepsilon})  $$
\item[(b)] for any sequence $\displaystyle \varepsilon_{n} \rightarrow 0$ and any 
 $\displaystyle (x,c)$, with $c: [0,1] \rightarrow U(x)$ 
$\displaystyle \mathring{\delta}^{x}$-derivable and 
$\displaystyle \mathring{Q}^{x}$-horizontal almost everywhere, 
there is a recovery sequence $\displaystyle (x_{n}, c_{n}) 
\in \mathcal{L}_{\varepsilon_{n}}(X,d,  \delta)$ such that 
$$l^{x}(c) \ = \ \lim_{n \rightarrow \infty} l^{x_{n}} (c_{n}) $$
\end{enumerate}

{\bf Proof of (a).} This is a consequence of propositions \ref{propcan}, 
\ref{phor} and definition \ref{defcoh} of a coherent projection. With the notations 
from (a) we see that we have to prove 
$$l^{x}(c) \ = \ \bar{l}^{x}(Q^{x} c)  \ \leq \ \liminf_{\varepsilon \rightarrow 0}
\bar{l}^{x_{\varepsilon}}(Q^{x_{\varepsilon}}_{\varepsilon} 
c_{\varepsilon})  $$
This is true because $\displaystyle (X, \bar{d}, \bar{\delta})$ is a tempered  dilatation 
structure and because of condition (A). Indeed 
 from the fact that $\displaystyle (X, \bar{d}, \bar{\delta})$ is tempered and from
 (\ref{uu}) (which is a consequence of condition (A)) we deduce that 
 $\displaystyle Q_{\varepsilon}$ is uniformly continuous on compact 
sets in a uniform way: for any compact set $K \subset X$ there is are constants 
$L(K) > 0$ (from (A)) and $C>0$ (from the tempered condition) such that for any 
$\varepsilon \in (0,1]$, any $x \in K$ and any $u,v$ sufficiently close to $x$ we
have: 
$$\bar{d}\left( Q^{x}_{\varepsilon} u, Q^{x}_{\varepsilon} v \right) \, \leq \, 
C \, \left(\bar{\delta}^{x}_{\varepsilon} \bar{d} \right) \left(
Q^{x}_{\varepsilon} u, Q^{x}_{\varepsilon} v \right) 
\, \leq \, C \, L(K) \, \bar{d}(u,v)$$
The sequence   $\displaystyle Q^{x}_{\varepsilon}$ uniformly converges to 
 $\displaystyle Q^{x} $ as $\varepsilon$ goes to $0$, 
uniformly with respect to $x$ in compact sets. Therefore if $\displaystyle 
 (x_{\varepsilon}, c_{\varepsilon}) \in \mathcal{L}_{\varepsilon}(X,d,  \delta)$
 converges to $(x,c)$ then $(x_{\varepsilon}, Q^{x_{\varepsilon}}_{\varepsilon} 
c_{\varepsilon}) \in \mathcal{L}_{\varepsilon}(X,\bar{d},  \bar{\delta})$ converges
to $(x,Q^{x}c)$. Use now the fact that by corollary \ref{cortemp} $(X, \bar{d}, \bar{\delta})$ is a length dilatation 
structure. The proof is done.

{\bf Proof of (b).} We have to construct a recovery sequence. We are doing
this by discretization of $c: [0,L] \rightarrow U(x)$. Recall that $c$ is a
curve which is $\displaystyle \mathring{\delta}^{x}$-derivable a.e. and 
$\displaystyle \mathring{Q}^{x}$-horizontal, that is for
almost every $t \in [0,L]$ the limit  
$$u(t) \ = \ \lim_{\mu \rightarrow 0} \delta^{x}_{\mu^{-1}} \, \Delta^{x} 
(c(t), c(t+\mu))$$ 
exists and $\displaystyle Q^{x} \, u(t) \, = \, u(t)$. Moreover we may suppose
that for almost every $t$ we have $\displaystyle \bar{d}^{x}(x, u(t)) \leq 1$ and 
$\displaystyle \bar{l}^{x}(c) \leq L$. 

There are functions $\displaystyle \omega^{1}, \omega^{2}: (0,+\infty)
\rightarrow [0,+\infty)$ with $\displaystyle \lim_{\lambda \rightarrow 0}
\omega^{i}(\lambda) = 0$, with the following property: 
for any  $\lambda > 0$ sufficiently small there is  a division $\displaystyle A_{\lambda} \, = \, \left\{ 0 < t_{0} < ... 
< t_{P} < L \right\}$ such that 
\begin{equation}
\frac{\lambda}{2} \, \leq \, \min \left\{ \frac{t_{0}}{t_{1} - t_{0}},
\frac{L-t_{P}}{t_{P} - t_{P-1}}, t_{k} - t_{k-1} \mbox{ : } k = 1, ... , P
\right\} 
\label{tpr1}
\end{equation}
\begin{equation}
\lambda \, \geq \, \max \left\{ \frac{t_{0}}{t_{1} - t_{0}},
\frac{L-t_{P}}{t_{P} - t_{P-1}}, t_{k} - t_{k-1} \mbox{ : } k = 1, ... , P
\right\} 
\label{tpr2}
\end{equation}
and such that $\displaystyle u(t_{k})$ exists for any $k = 1, ... , P$ and 
\begin{equation}
\mathring{d}^{x}(c(0), c(t_{0})) \, \leq \, t_{0} \, \leq \, \lambda^{2}
\label{tpr3}
\end{equation}
\begin{equation}
\mathring{d}^{x}(c(L), c(t_{P})) \, \leq \, L - t_{P} \, \leq \, \lambda^{2}
\label{tpr4}
\end{equation}
\begin{equation}
\mathring{d}^{x}(u(t_{k-1}), \Delta^{x}(c(t_{k-1}), c(t_{k})) \, \leq \, 
\left(t_{k} - t_{k-1} \right) \, \omega^{1}(\lambda)
\label{tpr5}
\end{equation}
\begin{equation}
\mid \int_{0}^{L} \bar{d}^{x}(x, u(t)) \mbox{ d}t \, - \, \sum^{P-1}_{k=0} 
(t_{k+1} - t_{k}) \, \bar{d}^{x}(x, u(t_{k})) \mid  \,  \leq \,
\omega^{2}(\lambda)
\label{tpr6}
\end{equation}
Indeed (\ref{tpr3}), (\ref{tpr4}) are a consequence of the fact that $c$ is 
$\displaystyle \mathring{d}^{x}$-Lipschitz, \eqref{tpr5} is a consequence of 
Egorov theorem applied to 
$$f_{\mu}(t) \ = \ \delta^{x}_{\mu^{-1}} \, \Delta^{x}(c(t), c(t+\mu))$$
and \eqref{tpr6} comes from the definition of the integral
$$l(c) \ = \ \int_{0}^{L} \bar{d}^{x}(x, u(t)) \mbox{ d}t $$
For each $\lambda$ we shall choose $\varepsilon = \varepsilon(\lambda)$ and we
shall construct a curve $\displaystyle c_{\lambda}$ with the properties: 
\begin{enumerate}
\item[(i)] $\displaystyle (x,c_{\lambda}) \, \in
\mathcal{L}_{\varepsilon(\lambda)}(X, d, \delta)$ 
\item[(ii)] $\displaystyle \lim_{\lambda \rightarrow 0}
l^{x}_{\varepsilon(\lambda)}(c_{\lambda}) \, = \, l^{x}(c)$. 
\end{enumerate}
At almost every $t$ the point $u(t)$ represents the velocity of the curve $c$ 
seen as the the left translation of $\frac{\mathring{d}^{x}}{dt} c(t)$ by the 
group operation $\displaystyle \Sigma^{x}(\cdot, \cdot)$ 
to $x$ (which is the neutral element for the mentioned operation). The
derivative (with respect to $\displaystyle \mathring{\delta}^{x}$) of the curve 
$c$ at $t$ is 
$$y(t) \, = \, \Sigma^{x}(c(t), u(t))$$

Let us take $\varepsilon > 0$, arbitrary for the moment. 
We shall use the points of the division $\displaystyle A_{\lambda}$ and 
for any $k = 0, ...,P-1$ we shall define the point: 
\begin{equation}
y^{\varepsilon}_{k} \ = \ \hat{Q}^{x, c(t_{k})}_{\varepsilon} \,
\Sigma^{x}_{\varepsilon}(c(t_{k}), u(t_{k}))
\label{tpr7}
\end{equation}
Thus $\displaystyle y^{\varepsilon}_{k}$ is obtained as the "projection" 
by $\displaystyle  \hat{Q}^{x, c(t_{k})}_{\varepsilon}$ of the "approximate 
left translation" $\displaystyle \Sigma^{x}_{\varepsilon}(c(t_{k}), \cdot)$ 
by $\displaystyle c(t_{k})$ of the velocity $\displaystyle u(t_{k})$. 
Define also the point: 
$$y_{k} \, = \, \Sigma^{x}(c(t_{k}), u(t_{k}))$$
By construction we have: 
\begin{equation}
y^{\varepsilon}_{k} \ = \ \hat{Q}^{x, c(t_{k})}_{\varepsilon} \,
y^{\varepsilon}_{k}
\label{tpr8}
\end{equation}
and by computation we see that $\displaystyle y^{\varepsilon}_{k}$ can be
expressed as: 
\begin{equation}
y^{\varepsilon}_{k} \ = \ \delta_{\varepsilon^{-1}}^{x} \,
Q^{\delta^{x}_{\varepsilon} c(t_{k})} \,
\delta^{\delta^{x}_{\varepsilon} c(t_{k})}_{\varepsilon} \, 
u(t_{k}) \ = \ 
\label{tpr9}
\end{equation}
$$= \ \Sigma^{x}_{\varepsilon}(c(t_{k}), Q^{\delta^{x}_{\varepsilon}
c(t_{k})} \, u(t_{k})) \ = \  \delta_{\varepsilon^{-1}}^{x} \,
\bar{\delta}^{\delta^{x}_{\varepsilon} c(t_{k})}_{\varepsilon} \,
Q^{\delta^{x}_{\varepsilon} c(t_{k})} \, u(t_{k})$$
Let us define the curve 
\begin{equation}
c^{\varepsilon}_{k}(s) \ = \ \hat{\delta}^{x, c(t_{k})}_{\varepsilon, s} \, 
y^{\varepsilon}_{k} \quad , \quad s \in [0, t_{k+1} - t_{k}]
\label{tpr11}
\end{equation}
which is a $\displaystyle \hat{Q}^{x}_{\varepsilon}$-horizontal curve 
(by supplementary hypothesis (B)) which joins $\displaystyle c(t_{k})$ with 
the point 
\begin{equation}
z^{\varepsilon}_{k} \ = \ \hat{\delta}^{x, c(t_{k})}_{\varepsilon, t_{k+1} - t_{k}} \,
y^{\varepsilon}_{k} 
\label{tpr10}
\end{equation}
The point $\displaystyle z^{\varepsilon}_{k}$ is an approximation of the point 
$$z_{k} \, = \, \mathring{\delta}^{x, c(t_{k})}_{t_{k+1} - t_{k}} y_{k}$$ 
We shall also consider the curve 
\begin{equation}
c_{k}(s) \ = \ \mathring{\delta}^{x, c(t_{k})}_{s} \, 
y_{k} \quad , \quad s \in [0, t_{k+1} - t_{k}]
\label{tpr111}
\end{equation}

There is a short curve $\displaystyle g^{\varepsilon}_{k}$  
which joins $\displaystyle 
z^{\varepsilon}_{k}$ with $\displaystyle c(t_{k+1})$, according to condition 
(Cgen). Indeed, for $\varepsilon$ sufficiently small the points 
$\displaystyle \delta^{x}_{\varepsilon} \, z^{\varepsilon}_{k}$ and 
$\displaystyle \delta^{x}_{\varepsilon} \, c(t_{k+1})$ are sufficiently close. 

Finally, take $\displaystyle g^{\varepsilon}_{0}$ and $\displaystyle
g^{\varepsilon}_{P+1}$ "short curves" which join $c(0)$ with 
$\displaystyle c(t_{0})$ and $\displaystyle c(t_{P})$ with $c(L)$ respectively. 

Correspondingly, we can find short curves $\displaystyle g_{k}$ (in the geometry of 
the dilatation structure $\displaystyle 
(U(x), \mathring{d}^{x}, \mathring{\delta}^{x}, \mathring{Q}^{x})$) joining 
$\displaystyle z_{k}$ with $\displaystyle c(t_{k+1})$, which are the uniform limit 
of the short curves $\displaystyle g^{\varepsilon}_{k}$ as $\varepsilon \rightarrow
0$. Moreover this convergence is uniform with respect to $k$ (and $\lambda$).
Indeed,  these short curves are made by $N$ curves of the type 
$\displaystyle s \mapsto  \hat{\delta}^{x, u_{\varepsilon}}_{\varepsilon, s}
v_{\varepsilon}$, with $\displaystyle \hat{Q}^{x,u_{\varepsilon}} v_{\varepsilon}
\, = \, v_{\varepsilon}$. Also, the short curves  $\displaystyle g_{k}$ are made
respectively by $N$  curves of the type $\displaystyle s \mapsto  
\mathring{\delta}^{x, u}_{s}
v$, with $\displaystyle \mathring{Q}^{x,u} v
\, = \, v$. Therefore we have:  
$$\bar{d}(\mathring{\delta}^{x, u}_{s} v , 
\hat{\delta}^{x,
 u_{\varepsilon}}_{\varepsilon, s} y^{\varepsilon}_{k}) \, = \, $$
$$\, = \, \bar{d} ( \Sigma^{x}(u, \bar{\delta}^{x}_{s} \Delta^{x}(u,v)) , 
\Sigma^{x}_{\varepsilon}(u_{\varepsilon}, \bar{\delta}^{\delta^{x}_{\varepsilon}
u_{\varepsilon}}_{s} \Delta^{x}_{\varepsilon}(u_{\varepsilon}, v_{\varepsilon})))$$ 
By an induction argument on the respective ends of segments forming the short
curves, using the axioms of coherent projections, we get the result.

By concatenation of all these curves we get two new curves: 
$$c^{\varepsilon}_{\lambda} \ = \ g^{\varepsilon}_{0} \left(\prod^{P-1}_{k = 0}
 c^{\varepsilon}_{k} \, g^{\varepsilon}_{k} \right) \,
g^{\varepsilon}_{P+1} $$
$$c_{\lambda} \ = \ g_{0} \left(\prod^{P-1}_{k = 0}
 c_{k} \, g_{k} \right) \,
g_{P+1} $$
From the previous reasoning we get that as $\varepsilon \rightarrow 0$ the curve 
$\displaystyle c^{\varepsilon}_{\lambda}$ uniformly converges to $\displaystyle 
c_{\lambda}$, uniformly with respect to $\lambda$. 

By theorem \ref{tancarnot}, specifically from relation (\ref{need22}) and
considerations below, we notice  that  for any $\displaystyle u \, = \, Q^{x} u$ 
the length of the curve $\displaystyle s \mapsto \delta^{x}_{s} u$ is: 
$$l^{x}(s \in [0,a] \mapsto \delta^{x}_{s} u) \, = \, a \, \bar{d}^{x}(x,u)$$
From here and relations (\ref{tpr3}), (\ref{tpr4}), (\ref{tpr5}), (\ref{tpr6}) we 
get that 
\begin{equation}
l^{x}(c) \, = \, \lim_{\lambda \rightarrow 0} l^{x}(c_{\lambda})
\label{needlam}
\end{equation}
 
Condition (B) and the fact that $\displaystyle (X, \bar{d}, \bar{\delta})$ is 
tempered imply that there is a positive function $\displaystyle \omega^{3}(\varepsilon) =
 \mathcal{O}(\varepsilon)$ such that 
 \begin{equation}
 \mid l^{x}_{\varepsilon}(c^{\varepsilon}_{\lambda}) - 
 l^{x}(c_{\lambda}) \mid \, \leq \, \frac{\omega^{3}(\varepsilon)}{\lambda}
 \label{needeps}
 \end{equation}
 This is true because if $\displaystyle v \, \hat{Q}^{x, u}_{\varepsilon} v$ then 
 $\displaystyle \delta^{x}_{\varepsilon} v \, = \, Q^{\delta^{x}_{\varepsilon} u} 
 \delta^{x}_{\varepsilon}v$, therefore by condition 
 (B)
 $$\frac{l^{x}_{\varepsilon}(s \in [0,a] \mapsto \hat{\delta}^{x, u}_{\varepsilon,
 s} v)}{\delta^{x}_{\varepsilon} \bar{d} (u, v)} \, = \, 
\frac{\bar{l}(s \in [0,a] \mapsto \bar{\delta}^{\delta^{x}_{\varepsilon} u}_{s} 
\delta^{x}_{\varepsilon} v)}{ \bar{d} (\delta^{x}_{\varepsilon} u ,
\delta^{x}_{\varepsilon} v)} \, \leq \, \mathcal{O}(\varepsilon) + 1$$
 Since each short curve is made by $N$ segments and the division $\displaystyle 
 A_{\lambda}$ is made by $1/\lambda$ segments, the relation (\ref{needeps})
 follows. 
 
 We shall choose now $\varepsilon(\lambda)$ such that $\displaystyle
 \omega^{3}(\varepsilon(\lambda)) \leq \lambda^{2}$ and we define: 
 $$c_{\lambda} \, = \, c^{\varepsilon(\lambda)}_{\lambda}$$
 These curves satisfy the properties (i), (ii). Indeed (i) is satisfied by 
 construction and (ii) follows from the choice of $\varepsilon(\lambda)$, 
 uniform convergence of $\displaystyle c^{\varepsilon}_{\lambda}$  to $\displaystyle 
c_{\lambda}$, uniformly with respect to $\lambda$, and relations (\ref{needeps}), 
and (\ref{needlam}).
 \hfill $\square$

\section{Conclusion}
\label{seccon}

In our opinion, the fact that sub-riemannian geometry may be described by 
about 12 axioms, {\bf without using any a priori given differential structure}, 
is remarkable and it shows the power of the dilatation structures approach. 
A geometry is not a simple object, for example euclidean geometry needs twice this
number of axioms. It should be clear that renouncing to such a basic object as a
differential structure is payed by the introduction of a number of axioms which
might seem too high at the first view. It is not so high though; just for an 
example, the number of axioms for the euclidean geometry decreases dramatically 
once we use as basic objects the algebraic and topological structure of real
numbers (or real vector spaces). 

Let us go back to Gromov viewpoint that the only intrinsic object of a
sub-riemannian space is the Carnot-Carath\'eodory distance. One of the 
most striking features of a regular sub-riemannian space is that it has at any
point a metric tangent space, which algebraically is a Carnot group. This has
been proved several times, by using the CC distance and lots of informations
coming from the underlying differential structure of the manifold. Let us
compare this with the result of Siebert, which characterizes homogeneous Lie
groups as locally compact groups admitting a contracting and continuous
one-parameter group of automorphisms. Siebert result has not a metric character.

In the work presented in this paper we tried to argue that we need more than
only the CC distance in order to describe regular sub-riemannian manifolds, but
less than the underlying differential structure: we need only dilatation
structures. Dilatation structures bring forth the other intrinsic ingredient,
namely the dilatations, which are generalizations of Siebert' contracting 
group of automorphisms.

\end{document}